\documentclass[12pt,a4paper]{article}
\usepackage[dvips]{epsfig}
\usepackage{amsmath}
\usepackage{amssymb}
\usepackage{amsfonts}
\usepackage{amsthm}
\usepackage{amsbsy}
\usepackage{amsgen}
\usepackage{amscd}
\usepackage{amsopn}
\usepackage{amstext}
\usepackage{amsxtra}

\newtheorem{theorem}{Theorem}[section]
\newtheorem{proposition}[theorem]{Proposition}
\newtheorem{lemma}[theorem]{Lemma}
\newtheorem{corollary}[theorem]{Corollary}

\theoremstyle{definition}

\newcommand {\Z} {\mathbb{Z}}
\newcommand {\ZZ} {\mathcal{Z}}
\newcommand {\R} {\mathbb{R}}

\newcommand {\Q} {\mathbb{Q}}
\newcommand {\1} {\mathbf{1}}

\newcommand {\A} {\mathcal {A}}

\newcommand {\N} {\mathbb {N}}
\newcommand {\M} {\mathcal {M}}
\newcommand {\PP} {\mathbb P}

\begin{document}

\title{The distribution of lattice points in elliptic annuli}

\author{Igor Wigman \\ School of Mathematical Sciences
\\ Tel Aviv University Tel Aviv 69978, Israel \\
e-mail: igorv@post.tau.ac.il}

\maketitle

\begin{abstract}
We study the distribution of the number of lattice points lying in
thin elliptical annuli. It has been conjectured by Bleher and
Lebowitz that, if the width of the annuli tend to zero while their
area tends to infinity, then the distribution of this number,
normalized to have zero mean and unit variance, is Gaussian. This
has been proved by Hughes and Rudnick for circular annuli whose
width shrink to zero sufficiently slowly. We prove this conjecture
for ellipses whose aspect ratio is transcendental and strongly
Diophantine, also assuming the width shrinks  slowly to zero.
\end{abstract}

\section{Introduction}

Let $B$ be an open convex domain in the plane containing the
origin, with a smooth boundary, and which is strictly convex (the
curvature of the boundary never vanishes). Let
\begin{equation*}
N_{B} (t):= \#Z^2 \cap t B,
\end{equation*}
be the number of integral points in the  $t$-dilate of $B$. As is
well-known, as $t\rightarrow\infty$, $N_{B}(t)$ is approximated by
the area of $t B$, that is
\begin{equation}
\label{eq:cnt fnc asym} N_{B}(t) \sim A t^2,
\end{equation}
where $A$ is the area of $B$.

A classical problem is to bound the size of the remainder
\begin{equation*}
\Delta_{B} (t) := N_{B} (t) - A t^2.
\end{equation*}
A simple geometric argument gives
\begin{equation}
\label{eq:Gss bnd} \Delta_{B} (t) = O(t),
\end{equation}
that is a bound in terms of the length of the boundary. It is
known that $\Delta_B$ is much smaller than the classical bound, as
Sierpinski ~\cite{SR} proved
\begin{equation*}
\Delta_{B} (t) = O(t^{2/3}).
\end{equation*}
Since then the exponent $2/3$ in this estimate has been improved
due to the works by many different researchers (see ~\cite{HUX}).
It is conjectured that one could replace the exponent by
$1/2+\epsilon$ for every $\epsilon>0$.

A different problem  is to study the value distribution of the
normalized error term, namely, of
\begin{equation*}
F_{B} (t) := \frac{\Delta_{B} (t)}{\sqrt{t}} = \frac{N_B (t) - A
t^2}{\sqrt{t}}.
\end{equation*}
Heath-Brown ~\cite{HB} treats this problem for $B = B(0,1)$, the
unit circle, and shows that there exists a probability density
$p(x)$, such that for every bounded continuous function $g(x)$,
\begin{equation*}
\lim\limits_{T\rightarrow\infty} \frac{1}{T} \int\limits_{0}^{T}
g(F_{B(0,1)} (t))dt = \int\limits_{-\infty}^{\infty} g(x) p(x) dx.
\end{equation*}
Somewhat surprisingly, the $p(t)$ is not a Gaussian: it decays as
$x\rightarrow\infty$ roughly as $\exp (-x^4)$, and it can be
extended to an entire function on a complex plane. Bleher
~\cite{BL3} establishes an analogue to Heath-Brown's theorem for
general ovals.

Motivated in parts by questions coming from mathematical physics,
we will concentrate on counting lattice points on annuli, namely,
integer points in
\begin{equation*}
(t+\rho) B \setminus t B,
\end{equation*}
that is, we study the remainder term of
\begin{equation*}
N_{B} (t, \, \rho) := N_{B} (t+\rho) - N_{B} (t),
\end{equation*}
where $\rho=\rho(t)$ is the width of annulus, depending on the
inner radius $t$. The "expected" number of points is the area $A
(2t\rho+\rho^2)$ of the annulus. Thus the corresponding {\em
normalized} remainder term is:
\begin{equation*}
S_{B} (t,\,\rho ) := \frac{N_{B} (t+\rho) - N_{B} (t) -A
(2t\rho+\rho^2) }{\sqrt{t}}
\end{equation*}

The statistics of $S_{B} (t,\,\rho )$ vary depending to the size
of $\rho(t)$. Of   particular interest to us are the following
regimes:

\noindent (1) The {\em microscopic regime} $\rho t$ is constant.
It was conjectured by Berry and Tabor ~\cite{BT} that the
statistics of $N_B (t,\,\rho)$ are Poissonian. Eskin, Margulis and
Mozes \cite{EMM} proved that the pair correlation function (which
is roughly equivalent to the variance of $N_B (t,\,\rho)$), is
consistent with the Poisson-random model.

\noindent (2) The {\em "global"}, or {\em "macroscopic"}, {\em
regime} $\rho(t)\rightarrow\infty$ (but $\rho=o(t)$). In such a
case, Bleher and Lebowitz ~\cite{BL4} showed that for a wide class
of $B$'s, $S_B (t,\,\rho)$ has a limiting distribution with tails
which decay roughly as $\exp(-x^4)$.

\noindent (3) The {\em intermediate} or {\em "mesoscopic"}, {\em
regime} $\rho\rightarrow 0$ (but $\rho t \rightarrow\infty$). If
$B$ is the inside of a "generic" ellipse
\begin{equation*}
\Gamma = \bigg\{ (x_1,\, x_2 ):\: x_1^2+ \alpha^2 x_2^2 = 1
\bigg\},
\end{equation*}
with $\alpha$ is Diophantine, the variance of $S_B (t,\,\rho)$ was
computed in ~\cite{BL} to be asymptotic to
\begin{equation}
\label{eq:sigma comp rho} \sigma^2 := \frac{8 \pi}{\alpha} \cdot
\rho
\end{equation}
For the circle ($\alpha=1$), the value is $16\rho
\log{\frac{1}{\rho}}$.

Bleher and Lebowitz \cite{BL4} conjectured that $S_B (t,\,\rho) /
\sigma$ has a standard Gaussian distribution. In 2004 Hughes and
Rudnick ~\cite{HR}  established the Gaussian distribution for the
unit circle, provided that $\rho(t) \gg t^{-\delta}$ for every
$\delta
> 0$.

In this paper, we prove the Gaussian distribution for the
normalized remainder term of "generic" {\em elliptic} annuli: We
say that $\alpha$ is {\em strongly Diophantine}, if for every
 $n \geq 1$ there is some $K>0$, such that for integers $a_j$
with $\sum\limits_{j=0}^{n} a_j \alpha^j \ne 0$,
\begin{equation*}
\bigg| \sum\limits_{j=0}^{n} a_j \alpha^j \bigg| \gg_{n}
\frac{1}{\bigg( \max\limits_{0\le j \le n} |a_j| \bigg) ^ {K}}.
\end{equation*}
This holds for any algebraic $\alpha$, for $\alpha = e$, and
almost every real $\alpha$, see section \ref{subsec:high mom}. Our
principal result is:

\begin{theorem}
\label{thm:norm dist} Let $B=\{x^2+\alpha^2 y^2 \le 1 \}$ with
$\alpha$ transcendental and strongly Diophantine. Assume that
$\rho = \rho (T) \rightarrow 0$, but for every $\delta > 0$, $\rho
\gg T^{-\delta}$. Then for every interval $\A$,
\begin{equation}
\lim_{T\rightarrow \infty } meas  \bigg\{ t\in [T,\, 2T ] :\:
\frac{S_B (t,\, \rho ) } {\sigma } \in \A \bigg\}  =
\frac{1}{\sqrt{2\pi }} \int\limits_{\A} e ^ {-\frac{x^2}{2}} dx,
\end{equation}
where $\sigma$ is given by \eqref{eq:sigma comp rho}.
\end{theorem}
This proves the conjecture of Bleher and Lebowitz in this case.

\paragraph{Remarks:} 1. In the formulation of theorem \ref{thm:norm dist} we assume for
technical reasons, that $\rho$ is a function of $T$ and
independent of $t\in [T,\, 2T]$. However one may easily see that
since $\rho$ may not decay rapidly, one may refine the result for
$\rho = \rho(t)$.

2. We compute statistics of the remainder term when the radius is
around $T$. A natural choice is assuming that the radius is
uniformly distributed in the interval $[T,\, 2T]$.

Our  case offers some marked differences from that of standard
circular annuli treated in \cite{HR}. To explain these, we note
that there are two main steps in treating these distribution
problems: The first step is to compute the moments of a
\underline{smoothed} version of $S_{B}$,  defined in
section~\ref{sec:smthng}.  We will show in section \ref{sec:smth
dist} that the moments of the smooth counting function are
Gaussian and that will suffice for establishing a normal
distribution for the smooth version of our problem. The second
step (section \ref{sec:unsmth}) is to recover the distribution of
the original counting function $S_{B}$ by estimating the variance
of the difference between $S_{B}$ and its smooth version. The
proof of that invokes a truncated Poisson summation formula for
the number of points of a general lattice which lie in a disk,
stated and proved in section~\ref{sec:asym for N_Lambda}.

The passage from circular annuli to general elliptical annuli
gives rise to new problems in both steps. The reason is that to
study the counting functions one uses Poisson summation to express
the counting functions as a sum over a certain lattice, that is as
a sum over closed geodesics of the corresponding flat torus.
Unlike the case of the circle, for a generic ellipse the sum is
over a lattice where the squared lengths of vectors are no longer
integers but of the form $n^2+m^2 \alpha^{-2}$, where $n,\, m \in
\Z$ and $\alpha$ is the aspect ratio of the ellipse.

One new feature present in this case is that these lengths can
{\em cluster} together, or, more generally, one may approximate
zero too well by the means of linear combinations of lengths. This
causes difficulties both in bounding the variance between the
original counting function and its smoothed version, especially in
the truncated summation formula of section~\ref{sec:asym for
N_Lambda}, and in showing that the moments of the smooth counting
function are given by "diagonal-like" contributions. This
clustering can be controlled when $\alpha$ is strongly
Diophantine.

Another problem we have to face, in evaluating moments of the
smooth counting function, is the possibility of non-trivial
correlations in the length spectrum. Their possible existence
(e.g. in the case of algebraic aspect ratio) obscures the nature
of the main term (the diagonal-like contribution) at this time. If
$\alpha$ is transcendental this problem can be overcome, see
proposition \ref{prop:diag impl prnc diag}.

\section{Smoothing}\label{sec:smthng}
Rather than counting integral points inside
elliptic annuli, we will count $\Lambda$-points inside
$B(0,1)$-annuli, where $\Lambda$ is a {\em lattice}. Denote the
corresponding counting function $N_{\Lambda}$, that is,
\begin{equation*}
N_{\Lambda} = \#\{\vec{n}\in\Lambda :\: |\vec{n}|\le t \}.
\end{equation*}

Let $\Lambda = \langle 1, \, i\alpha \rangle$ be a rectangular
lattice with $\alpha>0$ {\em transcendental} and {\em strongly
Diophantine} real number (almost all real $\alpha$ satisfy this,
see section \ref{subsec:high mom}). Denote
\begin{equation}
S_{\Lambda} (t,\, \rho ) = \frac { N_{\Lambda } (t + \rho ) -
N_{\Lambda } (t) - \frac{\pi}{d} (2 t \rho + \rho^{2})} {\sqrt {t}
}
\end{equation}
with $d := \det(\Lambda) = \alpha$. Thus
\begin{equation*}
S_{\Lambda} (t,\, \rho ) = S_B (t,\, \rho )
\end{equation*}
for an ellipse $B$ as in theorem \ref{thm:norm dist}, and we will
prove the result for $S_{\Lambda} (t,\, \rho )$.

We apply the same smoothing as in ~\cite{HR}: let $\chi$ be the
indicator function of the unit disc and $\psi$ a nonnegative,
smooth, even function on the real line, of total mass unity, whose
Fourier transform, $\hat{\psi}$ is smooth and has compact support
\footnote{To construct such a function, just take a function
$\phi$ with compact support and set $\hat{\psi} = \phi * \phi^*$
where $\phi^*(y) := \overline{\phi(-y)}$. Then $\psi =
|\check{\phi}|^2$ is nonnegative.}. One should notice that
\begin{equation}
\label{eq:cnt char} N_{\Lambda} (t) = \sum\limits_{\vec{n} \in
\Lambda} \chi \bigg( \frac{\vec{n}}{t} \bigg).
\end{equation}
Introduce a rotationally symmetric function $\Psi$ on $\R^2$ by
setting $\hat{\Psi} (\vec{y}) = \hat{\psi}(|\vec{y} | )$, where
$|\cdot |$ denotes the standard Euclidian norm. For $\epsilon >
0$, set
\begin{equation*}
\Psi_{\epsilon } (\vec{x} ) = \frac{1}{\epsilon ^ 2} \Psi \bigg(
\frac{\vec{x}}{\epsilon }\bigg).
\end{equation*}
Define in analogy with \eqref{eq:cnt char} a {\em smooth} counting
function
\begin{equation}
\label{eq:cnt char smth} \tilde{N}_{\Lambda, M} (t) =
\sum\limits_{\vec{n} \in \Lambda } \chi_{\epsilon } (
\frac{\vec{n}} {t} ),
\end{equation}
with $\epsilon = \epsilon (M)$, $\chi_{\epsilon } = \chi *
\Psi_{\epsilon }$, the convolution of $\chi$ with
$\Psi_{\epsilon}$. In what will follow,
\begin{equation}
\label{eq:def eps} \epsilon = \frac{1}{t \sqrt {M}},
\end{equation}
where $M = M(T)$ is the smoothness parameter, which tends to
infinity with $t$.

We are interested in the distribution of
\begin{equation}
\label{eq:ann cnt smth} \tilde{S}_{\Lambda , \, M, \, L} (t) =
\frac { \tilde{N}_{\Lambda ,\, M} (t + \frac{1}{L} ) -
                                \tilde{N}_{\Lambda ,\, M}
                                (t) - \frac{\pi}{d} (\frac{2 t}{L} + \frac{1}{L^{2}})}{\sqrt {t} },
\end{equation}
which is the smooth version of $S_{\Lambda} (t,\, \rho )$. We
assume that for every $\delta > 0$, $L = L(T) = O(T^{\delta } )$,
which corresponds to the assumption of theorem \ref{thm:norm dist}
regarding $\rho := \frac{1}{L}$. However, we will work with a
smooth probability space rather than just the Lebesgue measure.
For this purpose, introduce $\omega \ge 0$, a smooth function of
total mass unity, \ such that both $\omega$ and $\hat{\omega}$ are
rapidly decaying, namely
\begin{equation*}
|\omega(t) | \ll \frac{1}{(1+|t|)^A} , \; \; | \hat{\omega}(t) |
\ll \frac{1}{(1+|t|)^A},
\end{equation*}
for every $A>0$.

Define the averaging operator
\begin{equation*}
\langle f \rangle _T = \frac{1}{T} \int\limits_{-\infty}^{\infty}
f(t) \omega (\frac{t}{T} ) dt ,
\end{equation*}
and let $\PP_{\omega , \, T }$ be the associated probability
measure:
\begin{equation*}
\PP_{\omega , \, T } (f \in \A ) = \frac{1}{T}
\int\limits_{-\infty}^{\infty} 1_{\A} (f(t)) \omega (\frac{t}{T} )
dt ,
\end{equation*}

We will prove the following theorem in section \ref{sec:smth
dist}.

\begin{theorem}
\label{thm:norm dist smth} Suppose that $M(T)$ and $L(T)$ are
increasing to infinity with $T$, such that $M=O(T^{\delta})$ for
all $\delta > 0$, and $L/\sqrt{M} \rightarrow 0$. Then if $\alpha$
is {\em transcendental} and {\em strongly Diophantine}, we have
for $\Lambda = <1, \, i \alpha >$,
\begin{equation*}
\lim_{T\rightarrow \infty } \PP_{\omega,\, T}  \bigg\{
\frac{\tilde{S}_{\Lambda, \, M, \, L} } {\sigma } \in \A \bigg\}
= \frac{1}{\sqrt{2\pi }} \int\limits_{\A} e ^ {-\frac{x^2}{2}} d
\end{equation*}
for any interval $\A$, where $\sigma^2 := \frac{8 \pi }{d L} $.

\end{theorem}

\section{The distribution of $\tilde{S}_{\Lambda , \, M , \, L}$ }
\label{sec:smth dist}

We start from a well-known definition.
\paragraph{Definition:} A number $\mu$ is called {\em Diophantine}, if $\exists K > 0$, such that for
a rational $p/q$,
\begin{equation}
\label{eq:dioph num} \bigg| \mu - \frac{p}{q} \bigg| \gg_{\mu}
\frac{1}{q^{K}},
\end{equation}
where the constant involved in the $"\gg"$-notation depends only
on $\mu$. Khintchine proved that {\em almost all} real numbers are
Diophantine (see, e.g. ~\cite{S}, pages 60-63).

It is obvious from the definition, that $\mu$ is Diophantine iff
$\frac{1}{\mu}$ is such. For the rest of this section, we will
assume that $\Lambda ^{*} = \big\langle 1,\, i\beta \big\rangle$
with a Diophantine $\kappa := \beta ^ 2$, which satisfies
\eqref{eq:dioph num} with
\begin{equation}
\label{eq:k=k0} K=K_0,
\end{equation}
where $\Lambda^{*}$ is the {\em dual} lattice, that is $\beta :=
\frac{1}{\alpha}$. We may assume the Diophantinity of $\kappa$,
since theorem \ref{thm:norm dist} (and theorem \ref{thm:norm dist
smth}) assume $\alpha$'s being {\em strongly Diophantine}, which
implies, in particular, Diophantinity of $\alpha,\beta$ and
$\kappa$ (see the definition later in this section).

We will need a generalization of lemma 3.1 in ~\cite{HR} to a
general lattice $\Lambda$ rather than $\Z^2$.
\begin{lemma}
As $t\rightarrow \infty$,
\begin{equation}
\label{eq:app smth cnt} \tilde{N}_{\Lambda, M} (t) = \frac{\pi
t^2}{d} - \frac{\sqrt {t}}{d \pi} \sum\limits_{\vec{k} \in \Lambda
^ {*} \setminus \{ 0 \} } \frac{\cos \big( 2\pi t |\vec{k} | +
\frac{\pi}{4} \big) } {|\vec{k}|^{\frac{3}{2}}} \cdot
        \hat{\psi} \bigg( \frac{|\vec{k}|}{\sqrt{M}} \bigg) + O\bigg(\frac{1}{\sqrt {t}} \bigg),
\end{equation}
where, again $\Lambda^*$ is the dual lattice.
\end{lemma}

\begin{proof}
The proof is essentially the same as the one which obtains the
original lemma (see ~\cite{HR}, page 642). Using {\em Poisson
summation formula } on \eqref{eq:cnt char smth} and estimating
$\hat{\chi} (t \vec{k} )$ by the well-known asymptotics of the
Bessel $J_{1}$ function, we get:
\begin{equation*}
\tilde{N}_{\Lambda, M} (t) = \frac{\pi t^2}{d} - \frac{\sqrt
{t}}{d \pi} \sum\limits_{\vec{k} \in \Lambda ^ {*} \setminus \{ 0
\} } \bigg\{ \frac{\cos \big( 2\pi t |\vec{k} | + \frac{\pi}{4}
\big) } {|\vec{k}|^{\frac{3}{2}}} \cdot \hat{\psi} \bigg( \epsilon
t |\vec{k}| \bigg) + O \bigg( \frac{\hat{\psi} (\epsilon t
|\vec{k}| )} {t |\vec{k}|^{\frac{5}{2}}} \bigg) \bigg\},
\end{equation*}
where we get the main term for $\vec{k} = 0$. Finally, we obtain
\eqref{eq:app smth cnt} using \eqref{eq:def eps}. The contribution
of the error term is obtained due to the convergence of
$\sum\limits_{\vec{k} \in \Lambda ^ {*} \setminus \{ 0 \} }
\frac{1}{|\vec{k}|^{\frac{5}{2}}}$ as well as the fact that
$\hat{\psi} (x) \ll 1$.
\end{proof}

Unlike the standard lattice, if $\Lambda = \langle 1, \, i\alpha
\rangle$ with an irrational $\alpha^2$, then clearly there are no
nontrivial multiplicities, that is
\begin{lemma}
\label{lem:no mult} Let $\vec{a_i} = (n_i, \, m_i \cdot \alpha)
\in \Lambda$, $i = 1,\, 2$, with an irrational $\alpha^2$. If
$|\vec{a_1}| = |\vec{a_2}|$, then $n_1 = \pm n_2$ and $m_1 = \pm
m_2$.
\end{lemma}

By the definition of $\tilde{S}_{\Lambda, \, M ,\, L }$ in
\eqref{eq:ann cnt smth} and appropriately manipulating the sum in
\eqref{eq:app smth cnt} we obtain the following
\begin{corollary}
\begin{equation}
\begin{split}
\label{eq:approx S smth} \tilde{S}_{\Lambda, \, M ,\, L } (t) &=
\frac{2}{d \pi} \sum\limits_{\vec{k} \in \Lambda ^ {*} \setminus
\{ 0 \} } \frac{\sin \bigg( \frac{\pi |\vec{k}|}{L} \bigg) }
{|\vec{k}| ^{\frac{3}{2}}}
 \sin \bigg( 2 \pi \big( t+\frac{1}{2 L} \big) |\vec{k} | + \frac{\pi}{4} \bigg)
\hat{\psi} \bigg( \frac{|\vec{k}|}{\sqrt {M}} \bigg)  \\
&+ O \bigg(\frac{1}{\sqrt{t}} \bigg),
\end{split}
\end{equation}
\end{corollary}

We used
\begin{equation}
\label{eq:tailor sqrt} \sqrt{t+\frac{1}{L}} =\sqrt{t}
+O(\frac{1}{\sqrt{t}L})
\end{equation}
in order to change $\sqrt{t+\frac{1}{L}}$ multiplying the sum in
\eqref{eq:app smth cnt} for $N_{\Lambda } (t + \frac{1}{L} )$ by
$\sqrt{t}$. We use a smooth analogue of the simplest bound
\eqref{eq:Gss bnd} in order to bound the cost of this change to
the error term.

One should note that $\hat{\psi}$'s being compactly supported
means that the sum essentially truncates at $|\vec{k}| \approx
\sqrt{M}$.

\begin{proof}[Proof of theorem \ref{thm:norm dist smth}]

We will show that the moments of $\tilde{S}_{\Lambda, \, M ,\, L
}$ corresponding to the smooth probability space (i.e. $\langle
\tilde{S}_{\Lambda,\,M,\, L}^m \rangle_T $, see section
\ref{sec:smthng}) converge to the moments of the normal
distribution with zero mean and variance which is given by theorem
\ref{thm:norm dist smth}. This allows us to deduce that the
distribution of $\tilde{S}_{\Lambda, \, M ,\, L }$ converges to
the normal distribution as $T$ approaches infinity, precisely in
the sense of theorem \ref{thm:norm dist smth}.

First, we show that the mean is $O(\frac{1}{\sqrt{T}})$,
regardless of the Diophantine properties of $\alpha$. Since
$\omega$ is real,
\begin{equation*}
\Bigg| \Bigg\langle \sin \bigg( 2 \pi \big( t+\frac{1}{2 L} \big)
|\vec{k} | + \frac{\pi}{4} \bigg) \Bigg\rangle_T \Bigg| = \bigg|
\Im m \bigg\{ \hat{\omega} \big (-T |\vec{k}| \big) e ^ { i \pi
(\frac{|\vec{k}|}{L} + \frac{1}{4} } \bigg\} \bigg| \ll
\frac{1}{T^A |\vec{k}| ^ {A}}
\end{equation*}
for any $A>0$, where we have used the rapid decay of
$\hat{\omega}$. Thus
\begin{equation*}
\bigg| \bigg\langle \tilde{S}_{\Lambda, \, M ,\, L }
\bigg\rangle_T \bigg| \ll \sum\limits_{\vec{k} \in \Lambda ^ {*}
\setminus \{ 0 \} } \frac{1}{T^A |\vec{k}|^ {A+3/2}} + O \bigg(
\frac{1}{\sqrt{T}} \bigg) \ll  O \bigg( \frac{1}{\sqrt{T}} \bigg),
\end{equation*}
due to the convergence of $\sum\limits_{\vec{k} \in \Lambda ^ {*}
\setminus \{ 0 \} } \frac{1}{|\vec{k}|^ {A+3/2}}$, for $A >
\frac{1}{2}$

Now define
\begin{equation}
\label{eq:def M} \M_{\Lambda, \, m} := \Bigg\langle \bigg(
\frac{2}{d \pi} \sum\limits_{\vec{k} \in \Lambda ^ {*} \setminus
\{ 0 \} } \frac{\sin \bigg( \frac{\pi |\vec{k}|}{L} \bigg) }
{|\vec{k}| ^{\frac{3}{2}}}
 \sin \bigg( 2 \pi \big( t+\frac{1}{2 L} \big) |\vec{k} | + \frac{\pi}{4} \bigg)
\hat{\psi} \big( \frac{|\vec{k}|}{\sqrt {M}} \big) \bigg) ^{m}
\Bigg\rangle_T
\end{equation}

Then from \eqref{eq:approx S smth}, the binomial formula and the
Cauchy-Schwartz inequality,
\begin{equation*}
\bigg\langle \big( \tilde{S}_{\Lambda, \, M ,\, L }  \big) ^m
\bigg\rangle_T = \M_{\Lambda, \, m} + O \bigg( \sum
\limits_{j=1}^{m} \binom{m}{j} \frac{\sqrt {\M_{2m-2j}} }{T^{j/2}}
\bigg)
\end{equation*}

Proposition \ref{prop:var comp} together with proposition
\ref{prop:high mom comp} allow us to deduce the result of theorem
\ref{thm:norm dist smth} for a {\em transcendental strongly
Diophantine} $\beta^2$. Clearly, $\alpha$'s being transcendental
strongly Diophantine is sufficient.
\end{proof}

\subsection{The variance}
\label{ssec:var comp} The variance was first computed by Bleher
and Lebowitz ~\cite{BL} and we will give a version suitable for
our purpose. This will help the reader to understand our
computation of higher moments.
\begin{proposition}
\label{prop:var comp} Let $\alpha$ be Diophantine and $\Lambda =
\langle 1, \, i\alpha \rangle$. Then if for some fixed $\delta
> 0$, $M = O \big( T^{\frac{1}{K_0 + 1/2 + \delta}} \big)$ as $T\rightarrow\infty$,
then
\begin{equation*}
\bigg\langle \big( \tilde{S}_{\Lambda, \, M ,\, L }  \big) ^2
\bigg\rangle_T \sim  \sigma^2 :=
 \frac{2}{d^2 \pi^2}
\sum\limits_{\vec{k} \in \Lambda ^ {*} \setminus \{ 0 \} }
r(\vec{k})\frac{\sin ^2 \bigg( \frac{\pi |\vec{k}|}{L} \bigg) }
{|\vec{k}| ^{3}} \hat{\psi} ^ 2 \bigg( \frac{|\vec{k}|}{\sqrt {M}}
\bigg),
\end{equation*}
where
\begin{equation}
\label{eq:mult def} r(\vec{n}) =
\begin{cases}
1, \; &\vec{n} = (0,\,0) \\
2, \; &\vec{n} = (x,\, 0) \text{ or } (0,\, y) \\
4, &otherwise
\end{cases},
\end{equation}
is the ''multiplicity'' of $|\vec{n}|$. Moreover, if $L
\rightarrow \infty$, but $L/ \sqrt {M} \rightarrow 0$, then
\begin{equation}
\label{eq:sigma comp} \sigma^2 \sim \frac{8 \pi }{d L}
\end{equation}
\end{proposition}

\begin{proof}
Expanding out \eqref{eq:def M}, we have
\begin{equation}
\label{eq:M exp}
\begin{split}
\M_{\Lambda, \, 2} :=  \frac{4}{d^2 \pi ^2} &
\sum\limits_{\vec{k},\vec{l} \in \Lambda ^ {*} \setminus \{ 0 \} }
 \frac{\sin \bigg( \frac{\pi |\vec{k}|}{L} \bigg)
\sin \bigg( \frac{\pi |\vec{l}|}{L} \bigg) \hat{\psi} \big(
\frac{|\vec{k}|}{\sqrt {M}} \big) \hat{\psi} \big(
\frac{|\vec{l}|}{\sqrt {M}} \big) }
{|\vec{k}| ^{\frac{3}{2}} |\vec{l}| ^{\frac{3}{2}} } \\
& \times \bigg\langle \sin \bigg( 2 \pi \bigg( t+\frac{1}{2 L}
\bigg) |\vec{k} | + \frac{\pi}{4} \bigg) \sin \bigg( 2 \pi \bigg(
t+\frac{1}{2 L} \bigg) |\vec{l} | + \frac{\pi}{4} \bigg)
\bigg\rangle_T
\end{split}
\end{equation}
Now, it is easy to check that the average of the second line of
the previous equation is:
\begin{equation}
\label{eq:mean sines}
\begin{split}
\frac{1}{4}\bigg[
&\hat{\omega} \big( T(|\vec{k}| - |\vec{l}|) \big) e^{i\pi (1/L) (|\vec{l}| - |\vec{k}|)} + \\
&\hat{\omega} \big( T(|\vec{l}| - |\vec{k}|) \big) e^{i\pi (1/L) (|\vec{k}| - |\vec{l}|)} + \\
&\hat{\omega} \big( T(|\vec{k}| + |\vec{l}|) \big) e^{-i\pi (1/2 + (1/L) (|\vec{k}| + |\vec{l}|))} - \\
&\hat{\omega} \big( -T(|\vec{k}| + |\vec{l}|) \big) e^{i\pi (1/2 +
(1/L) (|\vec{k}| + |\vec{l}|))} \bigg]
\end{split}
\end{equation}
Recall that the support condition on $\hat{\psi}$ means that
$\vec{k}$ and $\vec{l}$ are both constrained to be of length
$O(\sqrt{M} )$, and so the off-diagonal contribution (that is for
$|\vec{k}| \ne |\vec{l}|$ ) of the first two lines of
\eqref{eq:mean sines} is
\begin{equation*}
\ll \sum\limits_ {\substack {\vec{k},\vec{l} \in \Lambda ^ {*}
\setminus \{ 0 \} \\ |\vec{k}|,\,|\vec{k'}| \le \sqrt {M}}}
\frac{M^{A(K_0+1/2)}}{T^A} \ll \frac{M^{A(K_0+1/2) + 2}}{T^A} \ll
T^{-B},
\end{equation*}
for every $B>0$, using lemma \ref{prop:off-diag low bnd}, the fact
that $|\vec{k}|,\, |\vec{l}| \gg 1$, $\big|\hat{\psi} \big| \ll
1$, and the assumption regarding $M$. We may use lemma
\ref{prop:off-diag low bnd} since we have assumed in the beginning
of this section that $\kappa$ is Diophantine.

Obviously, the contribution to \eqref{eq:M exp} of the two last
lines of \eqref{eq:mean sines} is negligible both in the diagonal
and off-diagonal cases, and so we are to evaluate the diagonal
approximation of \eqref{eq:M exp}, changing the second line of
\eqref{eq:M exp} by $1/2$, since the first two lines of
\eqref{eq:mean sines} are $2$. That proves the first statement of
the proposition. To find the asymptotics, we take a big parameter
$Y=Y(T)>0$ (which is to be chosen later), and write:
\begin{equation*}
\begin{split}
\sum\limits_{\substack{\vec{k},\,\vec{k'} \in \Lambda ^ {*}
\setminus \{ 0 \}  \\ |\vec{k}|=|\vec{k'}|}} \frac{\sin ^2 \bigg(
\frac{\pi |\vec{k}|}{L} \bigg) } {|\vec{k}| ^{3}} \hat{\psi} ^ 2
\bigg( \frac{|\vec{k}|}{\sqrt {M}} \bigg) &= \sum\limits_{\vec{k}
\in \Lambda ^ {*} \setminus \{ 0 \} } r(\vec{k} ) \frac{\sin ^2
\bigg( \frac{\pi |\vec{k}|}{L} \bigg) } {|\vec{k}| ^{3}}
\hat{\psi} ^ 2 \bigg( \frac{|\vec{k}|}{\sqrt {M}} \bigg) \\ &=
\sum\limits_{\substack{\vec{k} \in \Lambda ^ {*} \setminus \{ 0 \}
\\ |\vec{k}|^2 \le Y}} + \sum\limits_{\substack{\vec{k} \in
\Lambda ^ {*} \setminus \{ 0 \}  \\ |\vec{k}|^2 > Y}} := I_1 +
I_2,
\end{split}
\end{equation*}

Now for $Y=o(M)$, $\,\hat{\psi} ^ 2 \big( \frac{|\vec{k}|}{\sqrt
{M}} \big) \sim 1$ within the constraints of $I_1$, and so
\begin{equation*}
I_1 \sim \sum\limits_ {\substack{\vec{k} \in \Lambda ^ {*}
\setminus \{ 0 \}  \\ |\vec{k}|^2 \le Y}} r(\vec{k}) \frac{\sin ^2
\bigg( \frac{\pi |\vec{k}|}{L} \bigg) } {|\vec{k}| ^{3}}.
\end{equation*}
Here we may substitute $r(\vec{k})=4$, since the contribution of
vectors of the form $(x,\, 0)$ and $(0,\, y)$ is
$O(\frac{1}{L^2})$: representing their contribution as a
1-dimensional Riemann sum.

The sum in
\begin{equation*}
4 \sum\limits_ {\substack{\vec{k} \in \Lambda ^ {*} \setminus \{ 0
\}  \\ |\vec{k}|^2 \le Y}} \frac{\sin ^2 \bigg( \frac{\pi
|\vec{k}|}{L} \bigg) } {|\vec{k}| ^{3}}  = \frac{4 }{L}
\sum\limits_{\substack{\vec{k} \in \Lambda ^ {*} \setminus \{ 0 \}
\\ |\vec{k}|^2 \le Y}} \frac{\sin ^2 \bigg( \frac{\pi
|\vec{k}|}{L} \bigg) } {\big( \frac{|\vec{k}|}{L} \big) ^{3} }
\frac{1}{L^2}.
\end{equation*}
is a 2-dimensional Riemann sum of the integral
\begin{equation*}
\iint\limits_{1 / L^2 \ll x^2+ \kappa y^2 \le Y/L^2} \frac{\sin ^2
\big( \pi  \sqrt{x^2+ \kappa y^2} \big) } {|x^2+ \kappa y^2|
^{3/2}} dx dy \sim \frac{2\pi}{\beta}
\int\limits_{\frac{1}{L}}^{\frac{\sqrt{Y}}{L}} \frac{\sin ^2 (\pi
r)}{r^2} dr \rightarrow d \pi^3,
\end{equation*}
provided that $Y/L^2 \rightarrow \infty$, since
$\int\limits_{0}^{\infty} \frac{\sin^2 (\pi r)}{r^2} dr =
\frac{\pi^2}{2}$. We have changed the coordinates to the usual
elliptic ones. And so,
\begin{equation*}
I_1 \sim \frac{4 d \pi ^3}{L}
\end{equation*}
Next we will bound $I_2$. Since $\hat{\psi} \ll 1$, we may use the
same change of variables to obtain:
\begin{equation*}
I_2 \ll \frac{1}{L} \iint\limits_{x^2+ \kappa y^2 \ge Y/L^2}
\frac{\sin ^2 \big( \pi  \sqrt{x^2+ \kappa y^2} \big) } {|x^2+
\kappa y^2| ^{3/2}} dxdy \ll \frac{1}{L} \int \limits_{\sqrt{Y}/L}
^{\infty} \frac{dr}{r^2} = o \bigg( \frac{1}{L} \bigg).
\end{equation*}
This concludes the proposition, provided we have managed to choose
$Y$ with $L^2 = o(Y)$ and $Y = o(M)$. Such a choice is possible by
the assumption of the proposition regarding $L$.

\end{proof}

\begin{lemma}
\label{prop:off-diag low bnd} Suppose that $\vec{k},\,\vec{k'}
\in\Lambda^{*}$ with $|\vec{k}|,\,|\vec{k'}| \le \sqrt {M}$. Then
if $|\vec{k}| \ne |\vec{k'}|$,
\begin{equation*}
\big| |\vec{k}| - |\vec{k'}| \big| \gg M^ {-(K_0 + 1/2)}
\end{equation*}
\end{lemma}
\begin{proof}
\begin{equation*}
\big| |\vec{k}| - |\vec{k'}| \big| = \frac{\big| |\vec{k}|^2 -
|\vec{k'}|^2 \big|}{ |\vec{k}| + |\vec{k'}| } \gg
\frac{M^{-K_0}}{\sqrt{M}} = M^ {-(K_0 + 1/2)},
\end{equation*}
by \eqref{eq:dioph num} and \eqref{eq:k=k0}.
\end{proof}

\subsection{The higher moments}
\label{subsec:high mom} In order to compute the higher moments we
will prove that the main contribution comes from the so-called
{\em diagonal} terms (to be explained later). In order to be able
to bound the contribution of the {\em off-diagonal} terms, we
restrain ourselves to "generic" numbers, which are given in the
following definition:

\paragraph{Definition:} We call a number $\eta$ strongly
Diophantine, if it satisfies the following property: for any fixed
$n$, there exists $K_1\in\N$ such that for an integral polynomial
$P(x) = \sum\limits_{i=0}^n a_i x^i \in \Z[x]$, with $P(\eta) \ne
0$ we have
\begin{equation*}
\big| P(\eta ) \big| \gg_{\eta,\,n} h(P)^{-K_1},
\end{equation*}
where $h(P)=\max\limits_{0\le i \le n} {|a_i|}$ is the height of
$P$.

The fact that the strongly Diophantine numbers are "generic"
follows from various classical papers, e.g. \cite{MAH}.

Obviously, strong Diophantinity implies Diophantinity. Just as in
the case of Diophantine numbers $\eta$ is strongly Diophantine,
iff $\frac{1}{\eta}$ is such. Moreover, if $\eta$ is strongly
Diophantine, then so is $\eta^2$. As a concrete example of a
transcendental strongly Diophantine number, the inequality proven
by Baker ~\cite{BK} implies that for any {\em rational} $r\ne 0$,
$\eta = e^{r}$ satisfies the desired property.

We would like to make some brief comments concerning the number
$K_1$, which appears in the definition of a strongly Diophantine
number, although the form presented is sufficient for all our
purposes.

Let $\eta$ be a {\em real} number. One defines $\theta _k (\eta)$
to be $\frac{1}{k}$ times the supremum of the real numbers
$\omega$, such that $|P(\eta)| < h(P)^{-\omega}$ for infinitely
many polynomials $P$ of degree $k$. Clearly,
\begin{equation*}
\theta _k (\eta) = \frac{1}{k} \inf \{ \omega:\: |P(\eta)| \gg
_{\omega,\, k} h^{-\omega},\, \deg{P} = k \}.
\end{equation*}
It is well known ~\cite{TCH}, that $\theta _k (\eta) \ge 1$ for
all {\em transcendental} $\eta$. In 1932, Mahler \cite{MAH} proved
that $\theta_k (\eta) \le 4$ for almost all real $\eta$, and that
allows us to take any $K_1 > 4n$. He conjectured that
\begin{equation*}
 \theta_k (\eta) \le 1
\end{equation*}
which was proved in 1964 by Sprind\^{z}uk ~\cite{SP1},
~\cite{SP2}, making it legitimate to choose any $K_1 > n$.

Sprind\^{z}uk's result is analogous to Khintchin's theorem which
states that almost no $k$-tuple in $\R^k$ is {\em very well
approximable} (see e.g. ~\cite{S}, theorem 3A), for submanifold
$M\subset\R^k$, defined by
\begin{equation*}
M = \{(x,\,x^2,\,\ldots,\, x^k) :\: x\in \R \}.
\end{equation*}
The proof of this conjecture has eventually let to development of
a new branch in approximation theory, usually referred to as
"Diophantine approximation with dependent quantities" or
"Diophantine approximation on manifolds". A number of quite
general results were proved for a manifold $M$, see e.g.
~\cite{KM}.

We prove the following simple lemma which will eventually allow us
to exploit the strong Diophantinity of the aspect ratio of the
ellipse.
\begin{lemma}
\label{lem:gen bnd eps} If $\eta > 0$ is strongly Diophantine,
then it satisfies the following property: for any fixed natural
$m$, there exists $K \in \N$, such that if $z_j = a_j^2+\eta b_j^2
\ll M$, and $\epsilon_j = \pm 1$ for $j=1,\ldots , m$, with
integral $a_j,\, b_j$ and if $\sum\limits_{j=1}^{m} \epsilon_j
\sqrt{z_j} \ne 0$, then
\begin{equation}
\label{eq:gen bnd res} \big| \sum\limits_{j=1}^{m} \epsilon_j
\sqrt{z_j} \big| \gg_{\eta,\,m} M ^ {-K}.
\end{equation}
\end{lemma}

\begin{proof}
Let $m$ be given. We prove that every number $\eta$ that satisfies
the property of the definition of a strongly Diophantine number
with $n=2^{m-1}$, satisfies the inequality \eqref{eq:gen bnd res}
for some K, which will depend on $K_1$.

Let us $\{ \sqrt{z_j} \}_{j=1}^m$ be given. Suppose first, that
there is no $\{ \delta_j \} _{j=1} ^m \in \{ \pm 1 \} ^ m$ with
$\sum\limits_{j=1}^{m} \delta_j \sqrt{z_j} = 0$. Let us consider
\begin{equation*}
Q = Q(z_1,\,\ldots,\, z_m) := \prod\limits_{\{ \delta_j \} _{j=1}
^m \in \{ \pm 1 \} ^ m }
                              \sum\limits_{j=1}^{m} \delta_j \sqrt{z_j} \ne 0.
\end{equation*}
Now $Q = R \big( \sqrt {z_1}, \,\ldots,\, \sqrt{z_m} \big)$, where
\begin{equation*}
R \big( x_1,\,\ldots,\, x_m \big) := \prod\limits_{\{ \delta_j \}
_{j=1} ^m \in \{ \pm 1 \} ^ m }
                              \sum\limits_{j=1}^{m} \delta_j x_j.
\end{equation*}
Obviously, $R$ is a polynomial with integral coefficients of
degree $2^m$ such that for each vector $\underline{\delta} =
(\delta_j) = (\pm 1)$, $R ( \delta_1 x_1,\,\ldots,\, \delta_m x_m
) = R ( x_1,\,\ldots,\, x_m )$, and thus $Q (z_1,\,\ldots,\, z_m)$
is an integral polynomial in $z_1,\,\ldots,\,z_m$ of degree
$2^{m-1}$. Therefore, $Q = P(\eta )$, where $P$ is a polynomial of
degree $2^{m-1}$, $P=\sum\limits_{j=0}^{2^{m-1}} c_i x^{i}$, with
$c_i\in\Z$, such that $c_i = P_i (a_1,\,\ldots,\, a_m,\,
b_1,\,\ldots,\, b_m )$, where $P_i$ are polynomials. Thus there
exists $K_2$, such that $c_i \ll M^{K_2}$, and so, by the
definition of strongly Diophantine numbers, $Q \gg_{\eta,\,m} M^{-
K_2 K_1}$. We conclude the proof of lemma \ref{lem:gen bnd eps} in
this case by
\begin{equation*}
\big| \sum\limits_{j=1}^{m} \epsilon_j \sqrt{z_j}  \big| =
\frac{\big| Q \big| } {\bigg| \prod\limits_{\{ \delta_j \} _{j=1}
^m \ne \{ \epsilon_j \} _{j=1} ^m}
                              \sum\limits_{j=1}^{m} \delta_j \sqrt{z_j} \bigg|}
                              \gg_{\eta,\, m}
M^{-(K_2 K_1 + (2^{m}-1)/2)},
\end{equation*}
and so, setting $K:=K_2 K_1 + \frac{(2^{m}-1)}{2}$, we obtain the
result of the current lemma in this case.

Next, suppose that
\begin{equation}
\label{eq:subsum van} \sum\limits_{i=1}^{m} \delta_j \sqrt{z_j} =
0
\end{equation}
for some (given) $\{ \delta_i \}_{j=1}^{m} \in \{ \pm 1 \} ^m$.
Denote $S:= \{ j:\: \epsilon_j = \delta_j \}$, $S' =
\{1,\,\ldots,\, m \} \setminus S$. One should notice that
\begin{equation}
\label{eq:S,S', not empty} \emptyset \subsetneqq S,\,S'
\subsetneqq \{ 1,\, \ldots,\, m\}.
\end{equation}
Writing \eqref{eq:subsum van} in the new notations, we obtain:
\begin{equation*}
\sum\limits_{j\in S} \epsilon_j \sqrt{z_i} - \sum\limits_{j\in S'}
\epsilon_j \sqrt{z_i} = 0,
\end{equation*}
Finally,
\begin{equation*}
0 \ne \big| \sum\limits_{j=1}^{m} \epsilon_j \sqrt{z_j} \big| = 2
\big| \sum\limits_{j\in S'} \epsilon_j \sqrt{z_j} \big|
\gg_{\eta,\, m} M^{-K}
\end{equation*}
for some $K$ by induction, due to \eqref{eq:S,S', not empty}.
\end{proof}

\begin{proposition}
\label{prop:high mom comp} Let $m\in\N$ be given. Suppose that
$\alpha^2$ is \underline{transcendental} and \underline{strongly
Diophantine} which satisfy the property of lemma \ref{lem:gen bnd
eps} for the given $m$, with $K=K_m$. Denote $\Lambda = \langle
1,\, i\alpha \rangle$. Then if $\M = O \big(
T^{\frac{1-\delta}{K_m}} \big)$ for some $\delta > 0$, and if $L
\rightarrow \infty$ such that $L / \sqrt{M} \rightarrow 0$, the
following holds:
\begin{equation*}
\frac{\M_{\Lambda, \, m} }{\sigma^m} =
\begin{cases}
\frac{m!}{2^{m/2} \big( \frac{m}{2} \big) !  } + O\big( \frac{\log L}{L} \big), \; &m \text{ is even} \\
O\big( \frac{\log L}{L} \big), \; &m \text{ is odd}
\end{cases}
\end{equation*}
\end{proposition}

\begin{proof}
Expanding out \eqref{eq:def M},  we have
\begin{equation}
\label{eq:exp M_m}
\begin{split}
\M_{\Lambda, \, m} =  \frac{2^m}{d^m \pi ^m}
 \sum\limits_{\vec{k_1},\ldots ,\,\vec{k_m} \in \Lambda ^ {*} \setminus \{ 0 \} }
&\prod\limits_{j=1}^{m}
 \frac{\sin \bigg( \frac{\pi |\vec{k_j}|}{L} \bigg)
\hat{\psi} \big( \frac{|\vec{k_j}|}{\sqrt {M}} \big)}
{|\vec{k_j}| ^{\frac{3}{2}}} \\
& \times \bigg\langle \prod\limits_{j=1}^{m} \sin \bigg( 2 \pi
\big( t+\frac{1}{2 L} \big) |\vec{k_1} | + \frac{\pi}{4} \bigg)
\bigg\rangle_T
\end{split}
\end{equation}

Now,
\begin{equation*}
\begin{split}
\bigg\langle \prod\limits_{j=1}^{m}
&\sin \bigg( 2 \pi \big( t+\frac{1}{2 L} \big) |\vec{k_1} | + \frac{\pi}{4} \bigg) \bigg\rangle_T \\
&= \sum \limits _{\epsilon_j = \pm 1} \frac{\prod\limits_{j=1}^{m}
\epsilon_j} {2^m i^m} \hat{\omega} \bigg( -T \sum\limits_{j=1}^{m}
\epsilon_j |\vec{k_j}| \bigg) e^{ \pi i \sum\limits_{j=1} ^{m}
\epsilon_j \big( (1/L) |\vec{k_j}| + 1/4 \big) }
\end{split}
\end{equation*}

We call a term of the summation in \eqref{eq:exp M_m} with
$\sum\limits_{j=1}^{m} \epsilon_j |\vec{k_j}| = 0$ {\em diagonal},
and {\em off-diagonal} otherwise. Due to lemma \ref{lem:gen bnd
eps}, the contribution of the {\em off-diagonal} terms is:

\begin{equation*}
\ll \sum\limits_{\vec{k_1},\ldots ,\,\vec{k_m} \in \Lambda ^ {*}
\setminus \{ 0 \} } \bigg( \frac{T}{M ^ {K_m}} \bigg) ^ {-A} \ll
M^m T^{-A\delta},
\end{equation*}
for every $A>0$, by the rapid decay of $\hat{\omega}$ and our
assumption regarding $M$.

Since $m$ is constant, this allows us to reduce the sum to the
{\em diagonal terms}. The following definition and corollary
\ref{cor:diag impl prn diag} will allow us to actually sum over
the diagonal terms, making use of $\alpha$'s being transcendental.

\paragraph{Definition:}
We say that a term corresponding to $\{ \vec{k_1},\ldots,\,
\vec{k_m} \}  \in \bigg( \Lambda ^ {*} \setminus \{ 0 \} \bigg)
^m$ and $\{ \epsilon_j \} \in \{ \pm 1 \} ^m$ is a {\em principal
diagonal} term if there is a partition $\{ 1,\ldots,\, m \} =
\bigsqcup\limits_{i=1} ^{l} S_i$, such that for each $1\le i \le
l$ there exists a primitive $\vec{n_i}  \in  \Lambda ^ {*}
\setminus \{ 0 \}$, with non-negative coordinates, that satisfies
the following property: for every $j\in S_i$, there exist $f_j \in
\Z$ with $|\vec{k_j}| = f_j |\vec{n_i}|$. Moreover, for each $1\le
i \le l$, $\sum\limits_{j\in S_i} \epsilon_j f_j = 0$.

Obviously, the principal diagonal is contained within the
diagonal. However, if $\alpha$ is {\em transcendental}, the
converse is also true. It is easily seen, given the following
proposition.

\begin{proposition}
\label{prop:diag impl prnc diag} Suppose that $\eta\in\R$ is a
\underline{transcendental} number. Let
\begin{equation*}
z_j = a_j^2+\eta b_j^2
\end{equation*}
such that $(a_j,\, b_j ) \in \Z_+ ^2$ are all different primitive
vectors, for $1\le j \le m$. Then $\{ \sqrt{z_j} \} _{j=1} ^{m} $
are linearly independent over $\Q$.
\end{proposition}

The last proposition is an analogue of a well-known theorem due to
Besicovitch ~\cite{BS} about incommensurability of square roots of
integers. A proof of a much more general statement may be found
e.g. in ~\cite{BL2} (see lemma 2.3 and the appendix).

Thus we have
\begin{corollary}
\label{cor:diag impl prn diag} Every \underline{diagonal} term is
a \underline{principle diagonal} term whenether $\alpha$ is
\underline{transcendendal}.
\end{corollary}

By corollary \ref{cor:diag impl prn diag}, summing over diagonal
terms is the same as summing over {\em principal} diagonal terms.
Thus:
\begin{equation}
\label{eq:M from D}
\begin{split}
\frac{\M_{\Lambda, \, m}}{\sigma^m} &\sim \sum\limits_{l=1}^{m}
\sum\limits_{\{1,\ldots,\, m \} = \bigsqcup\limits_{i=1}^{l} S_i}
\bigg( \frac{1}{\sigma ^ {|S_1|}} \sum\limits_{\vec{n}_{1} \in
\Lambda ^ {*} \setminus \{ 0 \}}
{}' D_{\vec{n}_{1}} (S_1) \bigg) \\
\times &\bigg( \frac{1}{\sigma ^ {|S_2|}} \sum\limits_{ \Lambda ^
{*} \setminus \{ 0 \} \ni \vec{n}_2 \ne \vec{n}_1} {} '
D_{(\vec{n}_{2})} (S_2) \bigg) \ldots \bigg( \frac{1}{\sigma ^
{|S_l|}} \sum\limits_{ \Lambda ^ {*} \setminus \{ 0 \} \ni
\vec{n}_{l} \ne \vec{n}_2,\ldots,\, \vec{n}_{l-1}}{} '
D_{\vec{n}_l} (S_l) \bigg),
\end{split}
\end{equation}
where the inner summations are over primitive 1st-quadrant vectors
of $\Lambda ^ {*} \setminus \{ 0 \}$, and

\begin{equation*}
D_{\vec{n}} (S) = \frac { r(\vec{n}) }{|\vec{n}|^{3|S| / 2}}
\sum\limits_{\substack{f_j \ge 1 \\ \epsilon_j = \pm 1 \\
\sum\limits_{j\in S}  \epsilon_j f_j = 0}} \prod\limits_{j\in S}
\frac{-i\epsilon_j}{d\pi f_j^{3/2}} \sin \bigg( \frac{\pi}{L} f_j
|\vec{n}| \bigg) \hat{\psi} \bigg( \frac{|\vec{n}|}{\sqrt{M}}
\bigg) e^{i\pi\epsilon_j /4},
\end{equation*}
with $r(\vec{n})$ given by \eqref{eq:mult def}.

Lemma \ref{lem:sum D prim} allows us to deduce that the
contribution to \eqref{eq:M from D} of a partition is
$O(\frac{\log (L)}{L} )$, unless $|S_i| = 2$ for every $i =
1,\ldots,\, l$. In the latter case the contribution is 1 by the
2nd case of the same lemma. This is impossible for an odd $m$, and
so, it finishes the proof of the current proposition in that case.
Otherwise, suppose $m$ is even. Then the number of partitions
$\{1,\ldots,\, m \} = \bigsqcup\limits_{i=1}^{l} S_i$ with $|S_i|
= 2$ for every $1\le i \le l$ is
\begin{equation*}
\begin{split}
\frac{1}{\big( \frac{m}{2} \big) !} \binom{m}{2}\binom{m-2}{2}
\cdot \ldots \cdot \binom{2}{2} &= \frac{1}{\big( \frac{m}{2}
\big)!} \, \frac{m!}{2! \, (m-2)!} \, \frac{(m-2)!}{2!\, (m-4)!}
\cdot \ldots \cdot \frac{2!}{2!} \\ &=
 \frac{m!}{2^{m/2}\big( \frac{m}{2} \big)! }
\end{split}
\end{equation*}

That concludes the proof of proposition \ref{prop:high mom comp}.

\end{proof}

\begin{lemma}
\label{lem:sum D prim} If $L \rightarrow \infty$ such that $L /
\sqrt{M} \rightarrow 0$, then
\begin{equation*}
\frac{1}{\sigma^m} \bigg| \sum\limits_{\vec{n} \in \Lambda ^ {*}
\setminus \{ 0 \}} {} ' D_{\vec{n}} (S) \bigg| =
\begin{cases}
0,\;\; &|S|=1 \\
1, &|S|=2 \\
O \big( \frac{\log L }{L} ), &|S|\ge 3
\end{cases}
\end{equation*}
where the ' in the summation means that it is over primitive
vectors $(a, \, b)$.
\end{lemma}

\begin{proof}
Without loss of generality, we may assume that $S = \{ 1,\,2,\,
\ldots,\, |S| \}$, and we assume that $k:=|S|\ge 3$. Now,
\begin{equation}
\label{eq:sum D bnd} \bigg| \sum \limits_{\vec{n} \in \Lambda ^
{*} \setminus \{ 0 \}} {} ' D_{\vec{n}} (S) \bigg| \ll \sum
\limits_{\vec{n} \in \Lambda ^ {*} \setminus \{ 0 \}} {}
\frac{1}{|\vec{n}|^{3k/2}} Q(|\vec{n}|),
\end{equation}
where
\begin{equation*}
Q(z) := \sum\limits_{\{ \epsilon_j \} \in \{ \pm 1 \} ^ k}
\sum\limits_ {\substack{f_j \ge 1 \\ \sum\limits_{j=1}^{k}
\epsilon_j f_j = 0 }} \prod\limits_{j=1}^{k} \frac{| \sin (
\frac{\pi}{L} f_j z ) |}{f_j ^ {3/2} }.
\end{equation*}
Note that $Q(z) \ll 1$ for all $z$. We would like to establish a
sharper result for $z \ll L$. In order to have
$\sum\limits_{j=1}^{k} \epsilon_j f_j = 0$, at least two of the
$\epsilon_j$ must have different signs, and so, with no loss of
generality, we may assume, $\epsilon_k = -1$ and $\epsilon_{k-1} =
+1$. We notice that the last sum is, in fact, a Riemann sum, and
so
\begin{equation*}
\begin{split}
Q(z) \ll \frac{L^{k-1}}{L^{3k/2}} \int\limits_{1/L}^{\infty}
\cdots &\int\limits_{1/L}^{\infty} dx_1 \cdots dx_{k-2}
\sum\limits_{\{ \epsilon_j \}_{j=1} ^{k-2} \in \{ \pm 1 \}^{k-2}}
\int \limits_
{\frac{1}{L}+\max (0, -\sum\limits_{j=1}^{k-2} \epsilon_j f_j ) } ^{\infty} dx_{k-1} \\
\times & \bigg( \prod\limits_{j=1}^{k-1} \frac{\big|\sin ( \pi x_j
z )\big|}{x_j ^{3/2}} \bigg) \frac{\bigg|\sin \bigg( \pi z \cdot
\big(x_{k-1} + \sum\limits_{j=1}^{k-2} \epsilon_j x_j \big) \bigg)
\bigg| } {\bigg( x_{k-1}+\sum\limits_{j=1}^{k-1} \epsilon_j x_j
\bigg) ^ {3/2}}
\end{split}
\end{equation*}
By changing variables $y_i = z \cdot x_i$ of the last integral, we
obtain:
\begin{equation*}
\begin{split}
Q(z) \ll \frac{z^{k/2+1}}{L^{k/2 + 1}} \int\limits_{1/L}^{\infty}
\cdots &\int\limits_{1/L}^{\infty} dy_1 \cdots dy_{k-2}
\sum\limits_{\{ \epsilon_j \}_{j=1} ^{k-2} \in \{ \pm 1 \}^{k-2}}
\int \limits_
{\frac{z}{L}+\max (0, -\sum\limits_{j=1}^{k-2} \epsilon_j f_j ) } ^{\infty} dy_{k-1} \\
\times & \bigg( \prod\limits_{j=1}^{k-1} \frac{\big|\sin ( \pi y_j
)\big|}{y_j ^{3/2}} \bigg) \frac{\bigg|\sin \bigg( \pi\cdot
\big(y_{k-1} + \sum\limits_{j=1}^{k-2} \epsilon_j y_j \big) \bigg)
\bigg| } {\bigg( y_{k-1}+\sum\limits_{j=1}^{k-1} \epsilon_j y_j
\bigg) ^ {3/2}},
\end{split}
\end{equation*}

and since the last multiple integral is bounded, we may conclude
that
\begin{equation*}
Q(z) \ll
\begin{cases}
\frac{z^{k/2+1}}{L^{k/2+1}},\; &z < L \\
1, & z \ge L
\end{cases}
\end{equation*}
Thus, by \eqref{eq:sum D bnd},
\begin{equation*}
\bigg| \sum \limits_{\vec{n} \in \Lambda ^ {*} \setminus \{ 0 \}}
{} ' D_{\vec{n}} (S) \bigg| \ll
\sum \limits_{\substack{\vec{n} \in \Lambda ^ {*} \setminus \{ 0 \} \\
|\vec{n}| \le L}} \frac{1}{|\vec{n}|^ {3k/2}} \cdot
\frac{|\vec{n}| ^ {k/2+1}}{L^{k/2+1}} +
\sum \limits_{\substack{\vec{n} \in \Lambda ^ {*} \setminus \{ 0 \} \\
|\vec{n}| > L}} \frac{1}{|\vec{n}|^ {3k/2}} =: S_1+S_2.
\end{equation*}

Now, considering $S_1$ and $S_2$ as Riemann sums, and computing
the corresponding integrals in the usual elliptic coordinates we
get:

\begin{equation*}
S_1 \ll \frac{1}{L^{k/2+1}}
\sum \limits_{\substack{\vec{n} \in \Lambda ^ {*} \setminus \{ 0 \} \\
|\vec{n}| \le L}} \frac{1}{|\vec{n}|^{k-1}}
 \ll
\frac{1}{L^{k/2+1}} \int \limits_{1}^{L} \frac{dr}{r^{k-2}} \ll
\frac{\log L}{L^{k/2+1}},
\end{equation*}
since $k \ge 3$.

Similarly,
\begin{equation*}
S_2 \ll \int\limits_{L}^{\infty} \frac{dr}{r^{3k/2-1}} \ll
\frac{1}{L^{3k/2-2}} \ll  \frac{1}{L^{k/2 +(k-2)}} \ll
\frac{1}{L^{k/2+1}},
\end{equation*}
again since $k \ge 3$.

And so, returning to the original statement of the lemma, if $k =
|S| \ge 3$,
\begin{equation*}
\frac{1}{\sigma^m} \bigg| \sum\limits_{\vec{n} \in \Lambda ^ {*}
\setminus \{ 0 \}} {} ' D_{\vec{n}} (S) \bigg| \ll L^{k/2} \bigg(
\frac{\log L}{L^{k/2+1}} \bigg) \ll \frac{\log L}{L},
\end{equation*}
by \eqref{eq:sigma comp}.

In the case $|S|=2$, by the definition of $D_{\vec{n}}$ and
$\sigma^2$, we see that
\begin{equation*}
 \sum\limits_{\vec{n} \in \Lambda ^ {*} \setminus \{ 0 \}} {} ' D_{\vec{n}} (S) = \sigma^2.
\end{equation*}
This completes the proof of the lemma.

\end{proof}

\section{An asymptotical formula for $N_{\Lambda}$}
\label{sec:asym for N_Lambda}

We need  an asymptotical formula for the {\em sharp} counting
function $N_{\Lambda}$. Unlike the case of the standard lattice,
$\Z^2$, in order to have a good control over the error terms we
should use some Diophantine properties of the lattice we are
working with. We adapt the following notations:

Let $\Lambda$ be a lattice and $t > 0$ a real variable. Denote the
set of squared norms of $\Lambda$ by
\begin{equation*}
SN_{\Lambda} = \{ |\vec{n}|^2:\: n\in\Lambda \}.
\end{equation*}
Suppose we have a function $\delta_\Lambda
:SN_{\Lambda}\rightarrow \R$, such that given $\vec{k} \in \Lambda
$, there are \underline{no} vectors $\vec{n} \in \Lambda$ with $0
< | |\vec{n}|^2 - |\vec{k}|^2 | < \delta_{\Lambda} ( |\vec{k} | ^2
)$. That is,
\begin{equation*}
\Lambda \cap \{ \vec{n}\in\Lambda :\: |\vec{k}|^2 -
\delta_{\Lambda} ( |\vec{k} | ^2 ) < | \vec{n} |^2 < |\vec{k}|^2 +
\delta_{\Lambda} ( |\vec{k} | ^2 ) \} = A_{|\vec{k}|},
\end{equation*}
where
\begin{equation*}
A_{y} := \{ \vec{n}\in\Lambda :\: |\vec{n}| = y\}.
\end{equation*}
Extend $\delta_{\Lambda}$ to $\R$ by defining $\delta_{\Lambda}
(x) := \delta_{\Lambda} (|\vec{k}|^2 )$, where $\vec{k} \in
\Lambda$ minimizes $|x-|\vec{k} |^2|$ (in the case there is any
ambiguity, that is if $x = \frac{|\vec{n_1}|^2+|\vec{n_2}|^2}{2}$
for vectors $\vec{n_1},\,\vec{n_2} \in \Lambda$ with consecutive
increasing norms, choose $\vec{k}:= \vec{n_1}$). We have the
following lemma:

\begin{lemma}
\label{lem:asym cnt shrp} For every $a>0,\,c>1$,
\begin{equation*}
\begin{split}
N_{\Lambda } (t) &= \frac{\pi }{d} t^2 - \frac{\sqrt{t}}{d \pi }
\sum\limits_{\substack{\vec{k} \in \Lambda ^ {*} \setminus \{ 0 \}
\\ |\vec{k} | \le \sqrt{N} }} \frac{\cos \big( 2\pi t |\vec{k} | +
\frac{\pi}{4} \big) } {|\vec{k}| ^ \frac{3}{2}} +
O(N^{a}) \\
&+ O \bigg(\frac{t^{2c-1}}{\sqrt{N}} \bigg) + O
\bigg(\frac{t}{\sqrt{N}}
\cdot \big( \log t + \log (\delta_{\Lambda} (t^2) \big) \bigg) \\
&+ O \bigg( \log{N} + \log (\delta_{\Lambda ^ *} (t^2) ) \bigg)
\end{split}
\end{equation*}
\end{lemma}
As a typical example of such a function, $\delta_{\Lambda }$, for
$\Lambda = \langle 1, \, i\alpha \rangle$, with a Diophantine
$\gamma := \alpha ^2$, we may choose $\delta_{\Lambda } (y) =
\frac{c}{y^{K_0}}$, where $c$ is a constant. In this example, if
$\Lambda\ni\vec{k} = (a,b)$, then by lemma \ref{lem:no mult},
$A_{|\vec{k}|} = (\pm a,\, \pm b )$, provided that $\gamma$ is
irrational.

Our ultimate goal in this section is to prove lemma \ref{lem:asym
cnt shrp}. However, it would be more convenient to work with
$x=t^2$, and by abuse of notations we will call the counting
function $N_{\Lambda }$. Moreover, we will redefine
\begin{equation*}
N_{\Lambda } (x) :=
\begin{cases}
\# \{ \vec{k} :\: |\vec{k} | ^ 2 \le x \}, \; &x \ne |\vec{k}| ^2 \text{ for every } \vec{k}\in \Lambda \\
\# \{ \vec{k} :\: |\vec{k} | ^ 2 < x \} + 2, \; &\text{ otherwise}
\end{cases}
\end{equation*}
(recall that every norm of a $\Lambda$-vector is of multiplicity
4). We are repeating the argument of Titchmarsh ~\cite{TMHB} that
establishes the corresponding result for the remainder of the
arithmetic function, which counts the number of different ways to
write $m$ as a multiplication of a fixed number of natural
numbers.

Let $\Lambda = \big\langle 1, \, i\alpha \big\rangle$. For $\gamma
:= \alpha^2$, introduce a function $\ZZ_{\gamma} (s)$ (this is a
special value of an Eisenstein series) where $s= \sigma + it$ is a
complex variable. For $\sigma
> 1$, $\ZZ_{\gamma} (s)$ is defined by the following converging
series:
\begin{equation}
\label{eq:def of Z_gamma} \ZZ_{\gamma} (s) := \frac{1}{4}
\sum\limits_{\vec{k} \in \Lambda \setminus {0}}
\frac{1}{|\vec{k}|^{2s}}.
\end{equation}

Then $\ZZ_{\gamma}$ has an analytic continuation to the whole
complex plane, except for a single pole at $s=1$, defined by the
formula
\begin{equation*}
\Gamma (s) \pi^{-s} \ZZ_{\gamma} (s) = \int\limits_{1}^{\infty}
x^{s-1} \psi_{\gamma} (x) dx + \frac{1}{\sqrt{\gamma}}
\int\limits_{1}^{\infty} x^{-s} \psi_{1 / \gamma} (x) dx -
\frac{s-\sqrt{\gamma}(s-1)}{4 \sqrt{\gamma} s (1-s)},
\end{equation*}
where
\begin{equation*}
\psi_{\gamma} (x) := \frac{1}{4}\sum\limits_{\vec{k} \in \Lambda
\setminus {0}} e^{-\pi |\vec{k}|^2 x}.
\end{equation*}
This enables us to compute the residue of $\ZZ_{\gamma}$ at $s=1$:
\begin{equation*}
Res(\ZZ_{\gamma},\, 1) = \frac{\pi}{4 \sqrt{\gamma}}.
\end{equation*}
Moreover, $\ZZ_{\gamma}$ satisfies the following functional
equation:
\begin{equation}
\label{eq:func eq for Z_gamma} \ZZ_{\gamma} (s) =
\frac{1}{\sqrt{\gamma}} \chi (s) \ZZ_{1/\gamma} (1-s),
\end{equation}
with
\begin{equation}
\chi (s) = \pi ^{2s-1} \frac{\Gamma (1-s)}{\Gamma (s)}.
\end{equation}
We will adapt the notation
\begin{equation*}
\chi_{\gamma} (s) := \frac{1}{\sqrt{\gamma}} \chi (s).
\end{equation*}

The connection between $N_{\Lambda}$ and $\ZZ_{\gamma}$ is given
in the following formula, which is satisfied for every $c>1$:
\begin{equation*}
\frac{1}{4} N_{\Lambda} (x) = \frac{1}{2\pi i}
\int\limits_{c-i\infty}^{c+i\infty} \ZZ_{\gamma} (s) \frac{x^s}{s}
ds,
\end{equation*}
To prove it, just write $\ZZ_{\gamma}$ explicitly as the
converging series, and use
\begin{equation*}
\frac{1}{2\pi i} \int\limits_{c-i\infty}^{c+i\infty} \frac{y^s}{s}
ds = \nu (y),
\end{equation*}
where
\begin{equation*}
\nu (y) :=
\begin{cases}
1,\; &y > 1 \\
\frac{1}{2} & y=1 \\
0; \; & 0 < y < 1
\end{cases},
\end{equation*}
see ~\cite{D}, lemma on page 105, for example. One should bear in
mind that the infinite integral above is not converging, and so we
consider it in the symmetrical sense (that is, $\lim\limits_{T
\rightarrow \infty} \int\limits_{c-iT}^{c+iT}$ ).

The following lemma will convert the infinite vertical integral in
the last equation into a finite one, accumulating the
corresponding error term. It will make use of the Diophantine
properties of $\gamma$.

\begin{lemma}
\label{lem:inf to fin int} In the notations of lemma \ref{lem:asym
cnt shrp}, for any constant $c>1$,
\begin{equation}
\label{eq:inf to fin int} \frac{1}{4} N_{\Lambda} (x) =
\frac{1}{2\pi i} \int\limits_{c-i T}^{c+i T} \ZZ_{\gamma} (s)
\frac{x^s}{s} ds + O \bigg( \frac{x^c}{T} \bigg) +
 O \bigg( \frac{x }{T} \big( \log {x} + \log {\delta_{\Lambda} (x) }\big) \bigg)
\end{equation}
as $x,\, T \rightarrow \infty$.
\end{lemma}
\begin{proof}
Lemma on page 105 of ~\cite{D} asserts moreover that for $y \ne 1$
\begin{equation}
\label{eq:fin int approx to nu} \frac{1}{2\pi i} \int\limits_{c-i
T}^{c+i T} \frac{y^s}{s} ds = \nu (y) + O\bigg( y^c \min \big(1,\,
\frac{1}{T |\log {y} |} \big) \bigg),
\end{equation}
whereas for $y=1$,
\begin{equation}
\label{eq:fin int approx to nu y=1} \frac{1}{2\pi i}
\int\limits_{c-i T}^{c+i T} \frac{ds}{s} = \frac{1}{2} + O \bigg(
\frac{1}{T} \bigg)
\end{equation}

Suppose first that $x \ne |\vec{k}| ^2 \text{ for every }
\vec{k}\in \Lambda$. Summing \eqref{eq:fin int approx to nu} for
$y = \frac{x}{|\vec{k}|^2}$, where $\vec{k} \in \Lambda\setminus
\{0 \}$ gives (dividing both sides by 4):
\begin{equation*}
\frac{1}{2\pi i} \int\limits_{c-i T}^{c+i T} \ZZ_{\gamma} (s)
\frac{x^s}{s} ds = \frac{1}{4} N_{\Lambda} (x) + O \Bigg( x^c \sum
\limits_{\vec{k} \in \Lambda\setminus \{0 \}} \frac{\min{ \bigg(
1,\, \frac{1}{T \log { \frac{x}{| \vec{k}| ^{2}} }} \bigg)
}}{|\vec{k}|^{2c}} \Bigg).
\end{equation*}

The contribution to the error term of the right hand side of the
last equality of $\vec{k}\in \Lambda$ with $|\vec{k}|^2 > 2x$ or
$|\vec{k}|^2 < \frac{1}{2} x$ is
\begin{equation*}
\ll \frac{x^c}{T} \sum \limits_{|\vec{k}| \ge 2x \text{ or }
|\vec{k}| \le \frac{1}{2} x} \frac{1}{|\vec{k}| ^{2c}} \le
\frac{x^c}{T} \ZZ{_\gamma} (c) \ll \frac{x^c}{T}.
\end{equation*}

For vectors $\vec{k_0} \in \Lambda$, which minimize $\big|
|\vec{k}|^2 - x \big|$ (in the case of ambiguity we choose
$\vec{k_0}$ the same way we did in lemma \ref{lem:asym cnt shrp}
while extending $\delta_{\Lambda}$), the corresponding
contribution is
\begin{equation*}
\frac{x^c}{|\vec{k}|^{2c}} \ll \frac{x^c}{x^c} = 1.
\end{equation*}

Finally, we bound the contribution of vectors $\vec{k} \in \Lambda
\setminus \{ 0 \}$ with $|\vec{k_0}|^2 < |\vec{k}|^2 < 2x$, and
similarly, of vectors with $\frac{1}{2} x < |\vec{k}|^2 <
|\vec{k_0}|^2$. Now, by the definition of $\delta_{\Lambda}$,
every such $\vec{k}$ satisfies:
\begin{equation*}
|\vec{k}|^2 \ge |\vec{k_0}|^2 + \delta_{\Lambda} (x) \ge x +
\frac{1}{2} \delta_{\Lambda} (x).
\end{equation*}
Moreover, $\log {\frac{|\vec{k} |^2}{x}} \gg \frac{|\vec{k} |^2 -
|\vec{k_0}|^2}{x}$, and so the contribution is:
\begin{equation*}
\begin{split}
&\ll \frac{x^c}{x^c T} x \sum\limits_{x + \frac{1}{2}
\delta_{\Lambda} (x) \le |\vec{k}|^2 < 2x} \frac{1}{|\vec{k} |^2 -
|\vec{k_0}|^2  } \ll
\frac{x}{T} \int\limits_{\sqrt{|\vec{k_0}|^2 + \delta_{\Lambda} (x)} }^{\sqrt{2x}} \frac{r}{r^2 - |\vec{k_0}|^2} dr \\
&= \frac{x}{2 T} \int\limits_{|\vec{k_0}|^2 + \delta_{\Lambda}
(x)}^{2x} \frac{du}{u-|\vec{k_0}|^2} \ll \frac{x}{T} \log {\big(
u-|\vec{k_0}|^2 \big) } \bigg| _{|\vec{k_0}|^2 + \delta_{\Lambda}
(x)}^{2x}  \\ &\ll \frac{x}{T} \big( \log{x} +
\log{\delta_{\Lambda} (x)} \big)
\end{split}
\end{equation*}

If $x = |\vec{k_0}|^2$ for some $\vec{k_0}\in \Lambda$, the proof
is the same except that we should invoke \eqref{eq:fin int approx
to nu y=1} rather than \eqref{eq:fin int approx to nu} for
$|\vec{k}| = |\vec{k_0}|$.

That concludes the proof of lemma \ref{lem:inf to fin int}.
\end{proof}

\begin{proof}[Proof of lemma \ref{lem:asym cnt shrp}]
We use lemma \ref{lem:inf to fin int} and would like to move the
contour of the integral in \eqref{eq:inf to fin int} from $\sigma
= c$, $-T \le t \le T$ left to $\sigma = -a$ for some $a>0$. Now,
for $\sigma \ge c$,
\begin{equation*}
\big| \ZZ_{\gamma} (s) \big| = O \big( 1 \big),
\end{equation*}
and by the functional equation \eqref{eq:func eq for Z_gamma} and
the Stirling approximation formula,
\begin{equation*}
\big| \ZZ_{\gamma} (s) \big| \ll t^{1+2a}
\end{equation*}
for $\sigma = -a$. Thus by the Phragm\'{e}n-Lindel\"{o}f argument
\begin{equation*}
\big| \ZZ_{\gamma} (s) \big| \ll t^{(1+2a)(c-\sigma ) / (a+c)}
\end{equation*}
in the rectangle $-a-iT$, $c-iT$, $c+iT$, $-a+iT$. Using this
bound, we obtain
\begin{equation*}
\bigg| \int\limits_{-a+iT}^{c+iT} Z_{\gamma} (s) \frac{x^s}{s} ds
\bigg| \ll \frac{T^{2a}}{x^{a}} + \frac{x^c}{T},
\end{equation*}
and so is $\big| \int\limits _{-a-iT}^{c-iT} \big| $. Collecting
the residues at $s=1$ with residue being the main term of the
asymptotics,
\begin{equation*}
Res \big( Z_{\gamma} (s) \frac{x^s}{s} ,\, 1 \big) = \frac{\pi}{4
\sqrt{\gamma}} x
\end{equation*}
and at $s=0$ with
\begin{equation*}
Res \big( Z_{\gamma} (s) \frac{x^s}{s} ,\, 0 \big) =  Z_{\gamma}
(0) = O \big( 1 \big),
\end{equation*}
 we get:
\begin{equation*}
\begin{split}
\Delta_{\Lambda} (x) &:= \frac{1}{4} N_{\Lambda} (x) -
\frac{\pi}{4 \sqrt{\gamma}} x =
\frac{1}{2\pi i} \int\limits_{-a-i T}^{-a+i T} \ZZ_{\gamma} (s) \frac{x^s}{s} ds \\
&+ O \bigg( \frac{x^c}{T} \bigg) +
 O \bigg( \frac{x }{T} \big( \log {x} + \log {\delta_{\Lambda} (x) }\big) \bigg) +O(1) +
O \big( \frac{T^{2a}}{x^a} \big) .
\end{split}
\end{equation*}
Denote the integral in the last equality by $I$ and let $\kappa:=
\frac{1}{\gamma}$. Using the functional equation of $\ZZ_{\gamma}$
\eqref{eq:func eq for Z_gamma} again, and using the definition of
$\ZZ_{\kappa}$ for $\sigma > 1$, \eqref{eq:def of Z_gamma}, we
get:
\begin{equation}
\label{eq:int main term Delta} I=\frac{1}{2\pi i}
\int\limits_{-a-i T}^{-a+i T} \chi_{\gamma} (s) \ZZ_{\kappa} (1-s)
\frac{x^s}{s} ds = \frac{1}{2\pi i} \sum\limits_{\vec{k} \in
\Lambda^{*}} {}' \int\limits_{-a-i T}^{-a+i T} \frac{\chi_{\gamma}
(s)}{|\vec{k} | ^ {2-2s}} \frac{x^s}{s}ds,
\end{equation}
where the $'$ means that the summation is over vectors in the 1st
quadrant. Put
\begin{equation}
\label{eq:def of N} \frac{T^2}{\pi^2 x} := N+\frac{1}{2}
\delta_{\Lambda^*} (N),
\end{equation}
where $N = |\vec{k_0}|^2$ for some $\vec{k_0} \in \Lambda^{*}$ and
consider separately vectors $\vec{k} \in \Lambda^{*}$ with
$|\vec{k}|^2 > N$ and ones with $|\vec{k}|^2 \le N$.

First we bound the contribution of vectors $\vec{k}\in
\Lambda^{*}$ with $|\vec{k}|^2 > N$. Write the integral in
\eqref{eq:int main term Delta} as $\int\limits_{-a-i T}^{-a+i T} =
\int\limits_{-a-i T}^{-a-i} + \int\limits_{-a-i}^{-a+i}  +
\int\limits_{-a+i}^{-a+i T}$. Then
\begin{equation*}
\bigg| \sum\limits_{\vec{k} \in \Lambda^{*}} {}'
\int\limits_{-a-i}^{-a+i} \frac{\chi_{\gamma} (s)}{|\vec{k} | ^
{2-2s}} \frac{x^s}{s}ds \bigg| \ll x^{-a}
\sum\limits_{\substack{\vec{k} \in \Lambda^{*} \\ |\vec{k}|^2 >
N}} {}' \frac{1}{|\vec{k}|^{2+2a}} \le x^{-a} \ZZ_{\kappa} (1+a)
\ll x^{-a}.
\end{equation*}
Now,
\begin{equation*}
\begin{split}
|J| = \bigg| \int\limits_{-a+i}^{-a+i T} \frac{\chi_{\gamma}
(s)}{|\vec{k} | ^ {2-2s}} \frac{x^s}{s}ds \bigg| &=
 \frac{x^{-a} \pi^{-2a-1}}{\sqrt{\gamma} |\vec{k}| ^{2+2a}} \bigg| \int\limits_{1}^{T} i \frac{\Gamma (1-s)}{\Gamma(s)}
\frac{\big( |\vec{k}|^2 x \big) ^{ti}}{ti}  \pi^{2ti} dt \bigg| \\
&\ll \frac{x^{-a}}{|\vec{k}| ^{2+2a}} \bigg| \int\limits_{1}^{T}
e^{iF(t)} \bigg( t^{2a}+O \big( t^{2a-1} \big) \bigg) dt \bigg|,
\end{split}
\end{equation*}
with
\begin{equation*}
F(t) = 2t \big(-\log{t} +\log{\pi} +1 \big) +t \log \big(
|\vec{k}|^2 x \big) = t\log{\frac{\pi^2 e^2 |\vec{k}|^2 x}{t^2}},
\end{equation*}
due to the Stirling approximation formula.

One should notice that the contribution of the error term in the
last bound is
\begin{equation*}
\ll \frac{T^{2a}}{x^a} \sum\limits_{\vec{k} \in \Lambda^{*}}{} '
\frac{1}{|\vec{k}|^{2+2a}} = \frac{T^{2a}}{x^a} \ZZ_{\kappa} (1+a)
\ll N^a.
\end{equation*}

We would like to invoke lemma 4.3 of ~\cite{TMHB} in order to
bound the integral above. For this purpose we compute the
derivative:
\begin{equation*}
F' (t) =  \log\bigg( \frac{|\vec{k}|^2 x \pi^2}{t^2} \bigg) \ge
\log {\bigg( \frac{|\vec{k}|^2}{N+\frac{1}{2} \delta_{\Lambda^*}
(N)} \bigg)},
\end{equation*}
by the definition of $N$, \eqref{eq:def of N}. Thus in the
notations of lemma 4.3 of ~\cite{TMHB},
\begin{equation*}
\frac{F'(t)}{G(t)} = \frac{\log \bigg( \frac{|\vec{k}|^2 x \pi^2}
{t^2} \bigg) }{t^{2a}}  \ge \frac{\log {\bigg(
\frac{|\vec{k}|^2}{N+\frac{1}{2} \delta_{\Lambda^*} (N)}
\bigg)}}{T^{2a}}.
\end{equation*}
We would also like to check that $\frac{G(t)}{F'(t)}$ is
monotonic. Differentiating that function and leaving only the
numerator, we get:
\begin{equation*}
-t^{2a-1} \big( 2a\log \frac{|\vec{k}|^2 x \pi^2}{t^2} + 2 \big) <
-2a \cdot t^{2a-1} \log \frac{|\vec{k}|^2 }{N+\frac{1}{2}
\delta_{\Lambda^*} (N)}  < 0,
\end{equation*}
since $|\vec{k}^2| > N$. Thus
\begin{equation*}
|J| \ll \frac{x^{-a}}{|\vec{k}|^{2+2a}} \frac{T^{2a}}{\log
\frac{|\vec{k}|^2}{N+\frac{1}{2} \delta_{\Lambda^*} (N)}},
\end{equation*}
getting the same bound for $\big| \int\limits_{-a-i T}^{-a-i}
\big|$, and therefore we are estimating
\begin{equation*}
\sum\limits_{\vec{k} \in \Lambda^{*}}{} '  \frac{T^{2a}}{\log
\frac{|\vec{k}|^2}{N+\frac{1}{2} \delta_{\Lambda^*} (N)}}.
\end{equation*}
For $|\vec{k}|^2 \ge 2N$, the contribution of the sum in
\eqref{eq:int main term Delta} is:
\begin{equation*}
\ll \frac{T^{2a}}{x^a} \sum\limits_{\substack{\vec{k} \in
\Lambda^{*} \\ |\vec{k}|^2 \ge 2N}} \frac{1}{|\vec{k}|^{2+2a}} \le
\frac{T^{2a}}{x^a} \ZZ_{\kappa} (1+a) \ll N^{a}
\end{equation*}
As for vectors $\vec{k}\in\Lambda^*$ with $N+\delta_{\kappa} (N)
\le |\vec{k}|^2 < 2N$,
\begin{equation*}
\log \frac{|\vec{k}|^2}{N+\frac{1}{2}\delta_{\kappa} (N)} \gg
\frac{|\vec{k}|^2 - N}{N},
\end{equation*}
which implies that the corresponding contribution to the sum in
\eqref{eq:int main term Delta} is:
\begin{equation*}
\begin{split}
\ll \frac{T^{2a}}{x^a N^{1+a}} \sum\limits_{\substack{\vec{k} \in \Lambda^{*} \\
N+\delta_{\kappa} (N) \le |\vec{k}|^2 < 2N}} {} '
\frac{N}{|\vec{k}|^2 - N} &\ll
\int\limits_{\sqrt{N+\delta_{\kappa} (N)}}^{\sqrt{2N}} \frac{r}{r^2 - N - \frac{1}{2}\delta_{\kappa} (N)} \\
&\ll \log{\big(\delta_{\kappa} (N) \big)} + \log{N}
\end{split}
\end{equation*}

The main term of I comes from $|\vec{k}|^2 \le N$. For such a
$\vec{k}$, we write
\begin{equation}
\label{eq:mov int to im ax} \int\limits_{-a-i T}^{-a+i T} =
\int\limits_{-i \infty}^{i \infty} - \Bigg( \int\limits_{i T}^{i
\infty} + \int\limits_{-i \infty}^{-i T} + \int\limits_{-i
T}^{-a-i T} + \int\limits_{-a+i T}^{i T} \Bigg),
\end{equation}
that is, we are moving the contour of the integration to the
imaginary axis.

Consider the first integral in the brackets. It is a constant
multiple of
\begin{equation*}
 \int\limits_{T}^{\infty} e^{iF(t)} dt \ll \frac{1}{\log \bigg( \frac{N+\frac{1}{2}\delta_{\kappa} (N) }{|\vec{k}|^2} \bigg)},
\end{equation*}
and so the contribution of the corresponding sum is
\begin{equation*}
\begin{split}
\ll \sum\limits_{\substack{\vec{k} \in \Lambda^{*} \\
|\vec{k}|^2 \le N}} \frac{1}{|\vec{k}|^2  \log \bigg(
\frac{N+\frac{1}{2}\delta_{\kappa} (N) }{|\vec{k}|^2} \bigg)} &\ll
N\int\limits_{1}^{\sqrt{N}} \frac{dr}{r\big(
N+\frac{1}{2}\delta_{\kappa} (N) - r^2 \big)} \\ &\ll \log{N} +
\log{\delta_{\kappa} (N) },
\end{split}
\end{equation*}
by lemma 4.2 of ~\cite{TMHB}, and similarly for the second
integral in the brackets in \eqref{eq:mov int to im ax}.

The last two give
\begin{equation*}
\begin{split}
\ll \sum\limits_{\substack{\vec{k} \in \Lambda^{*} \\ |\vec{k}|^2
\le N}} \frac{1}{|\vec{k}|^2} \int\limits_{-a}^{0} \bigg(
\frac{|\vec{k}|^2 x}{T^2} \bigg) ^{\sigma} d\sigma &\ll
\sum\limits_{\substack{\vec{k} \in \Lambda^{*} \\ |\vec{k}|^2 \le
N}} \frac{1}{|\vec{k}|^2}
\bigg( \frac{T^2}{|\vec{k}|^2 x} \bigg) ^a \\
&\ll \frac{T^{2a}}{x^a} \int\limits_{1}^{\sqrt{N}}
\frac{dr}{r^{2a+1}} \ll N^a.
\end{split}
\end{equation*}

Altogether we have now proved:
\begin{equation}
\label{eq:asym form delt}
\begin{split}
\Delta_{\Lambda} (x) &= \frac{1}{2\pi^2 d i}
\sum\limits_{\substack{\vec{k} \in \Lambda^{*} \\ |\vec{k}|^2 \le
N}} {}' \frac{1}{|\vec{k}|^2} \int\limits_{- i \infty}^{i \infty}
\pi^{2s} \frac{\Gamma (1-s)}{\Gamma(s)} \frac{\big( |\vec{k}|^2 x
\big) ^s}{s} ds
+O(N^{a}) \\
&+ O \bigg(\frac{x^{c-1/2}}{\sqrt{N}} \bigg) +
O \bigg(\frac{\sqrt{x}}{\sqrt{N}} \cdot \big( \log x + \log (\delta_{\Lambda} (x) ) \big) \bigg) \\
&+ O \big( \log{N} + \log (\delta_{\Lambda ^ *} (x) ) \big)
\end{split}
\end{equation}

Recall the integral $\int\limits_{- i \infty}^{i \infty}
\frac{\Gamma (1-s)}{\Gamma(s)}\frac{y ^s}{s} ds $ is a principal
value, that is $\lim\limits_{T\rightarrow\infty} \int_{-iT}^{iT}$.
We have
$$\lim_{T\rightarrow\infty} \int_{-iT}^{iT} \frac{\Gamma (1-s)}
{\Gamma(s)}\frac{y ^s}{s} ds = -\sqrt{y} J_{-1}(2\sqrt{y})$$ as
can be seen by shifting contours. Note that the analogous
Barnes-Mellin formula
$$J_\nu(x) = \frac{1}{2\pi i} \int \limits_{-i\infty}^{i\infty} \Gamma(-s)
[\Gamma(\nu + s + 1)]^{-1} (x/2)^{\nu+2s} ds$$ valid for
$Re(\nu)>0$ (see ~\cite{EA}, (36), page 83), which deals with
convergent integrals, is proved in this manner.

The well-known asymptotics of the bessel $J$-function,
\begin{equation*}
J_{-1} (y) = \sqrt{\frac{2}{\pi y}} \cos\big(y+\frac{\pi}{4}\big)
+ O(y^{-3/2})
\end{equation*} as $y\rightarrow\infty$, allow us to estimate
the integral involved in \eqref{eq:asym form delt} in terms of x
and $\vec{k}$. Collecting all the constants and the error terms,
we obtain the result of lemma \ref{lem:asym cnt shrp}.
\end{proof}

\section{Unsmoothing}
\label{sec:unsmth}
\begin{proposition}
\label{prop:bnd 2nd mom diff} Let a lattice $\Lambda = \big\langle
1,\, i\alpha \big\rangle$ with a Diophantine $\gamma := \alpha^2$
be given. Suppose that $L \rightarrow \infty$ as $T\rightarrow
\infty$ and choose $M$, such that $L/\sqrt{M} \rightarrow 0$, but
$M = O \big( T^\delta \big)$ for every $\delta > 0$ as $T
\rightarrow \infty$. Suppose furthermore, that $M = O (L^{s_0})$
for some (fixed) $s_0 > 0$. Then
\begin{equation*}
\Bigg\langle \bigg| S_{\Lambda} (t,\, \rho ) - \tilde{S}_{\Lambda,
\, M, \, L} (t) \bigg| ^ {2} \Bigg\rangle_T \ll \frac{1}{\sqrt{M}}
\end{equation*}
\end{proposition}

\begin{proof}
Since $\gamma$ is Diophantine, we may invoke lemma \ref{lem:asym
cnt shrp} with $\delta_{\Lambda } (y) = \frac{c_1}{y^{K_0}}$ and
$\delta_{\Lambda ^* } (y) = \frac{c_2}{y^{K_0}}$, where $c_1,\,
c_2$ are constants. Choosing $a=\delta'$ and $c=1+\delta'/2$ for
$\delta'>0$ arbitrarily small and using essentially the same
manipulation we used in order to obtain \eqref{eq:approx S smth},
and using \eqref{eq:tailor sqrt} again, we get the following
asymptotical formula:
\begin{equation}
\label{eq:approx S shrp} S_{\Lambda} (t,\, \rho ) = \frac{2}{d
\pi} \sum\limits_{\vec{k} \in \Lambda ^ {*} \setminus \{ 0 \} }
\frac{\sin \bigg( \frac{\pi |\vec{k}|}{L} \bigg) } {|\vec{k}|
^{\frac{3}{2}}}
 \sin \bigg( 2 \pi \bigg( t+\frac{1}{2 L} \bigg) |\vec{k} | + \frac{\pi}{4} \bigg) + R_{\Lambda} (N,\,t ) ,
\end{equation}
where
\begin{equation*}
|R_{\Lambda} (N,\,t )| \ll \frac{N^{\delta'}}{\sqrt{|t|}} +
\frac{|t|^{1/2+\delta'}}{\sqrt{N}} + \frac{1}{|t|^{1/2-\delta'}}.
\end{equation*}

Set $N = T^{3}$. Since M is small, the infinite sum in
\eqref{eq:approx S smth} is truncated before $n = T^{3}$. Thus
\eqref{eq:approx S smth} together with \eqref{eq:approx S shrp}
implies:
\begin{equation}
\label{eq:diff asmp sum}
\begin{split}
&S_{\Lambda} (t,\, \rho ) -  \tilde{S}_{\Lambda, \, M, \, L} (t) = \\
\frac{2}{d \pi} &\sum\limits_{\substack{\vec{k} \in \Lambda ^ {*}
\setminus \{ 0 \} \\ |\vec{k}| \le T^{3/2}}} \frac{\sin \bigg(
\frac{\pi |\vec{k}|}{L} \bigg) } {|\vec{k}| ^{\frac{3}{2}}}
 \sin \bigg( 2 \pi \bigg( t+\frac{1}{2 L} \bigg) |\vec{k} | + \frac{\pi}{4} \bigg)
\bigg(1 - \hat{\psi} \big( \frac{|\vec{k}|}{\sqrt {M}} \big) \bigg)  \\
&\;\; \; \;\;\;\;\;\; \quad +R_{\Lambda} (T^{3},\, t ).
\end{split}
\end{equation}

Let $P_{\Lambda} (N,\, t)$ denote the sum in \eqref{eq:diff asmp
sum}. Then the Cauchy-Schwartz inequality gives:
\begin{equation}
\label{eq:CS ineq for var}
\begin{split}
\Bigg\langle \bigg| S_{\Lambda} (t,\, \rho ) - \tilde{S}_{\Lambda,
\, M, \, L} (t) \bigg| ^ {2}  \Bigg\rangle_T
 =
&\big\langle P_{\Lambda} ^2 \big\rangle_T +
\big\langle R_{\Lambda} (N,\, t) ^2 \big\rangle_T +  \\
&O \bigg( \sqrt{ \big\langle P_{\Lambda} ^2 \big\rangle_T}
\sqrt{\big\langle R_{\Lambda} (N,\, t) ^2 \big\rangle_T } \bigg).
\end{split}
\end{equation}

Observe that for the chosen $N$,
\begin{equation*}
\big\langle R_{\Lambda} (N,\,t ) ^ 2 \big\rangle_T = O \big(
T^{-1+\delta'}  \big)
\end{equation*}
for arbitrary small $\delta' > 0$, since the above equality is
satisfied pointwise.

Next we would like to bound $\big\langle P_{\Lambda} ^2
\big\rangle_T$. Just as we did while computing the variance of the
smoothed variable, $\tilde{S}_{\Lambda, \, M, \, L}$, we divide
all the terms of the expanded sum into the diagonal terms and the
off-diagonal ones (see section \ref{ssec:var comp}). Namely,
\begin{equation}
\label{eq:asmp eqlty P}
\begin{split}
\big\langle  P_{\Lambda} ^2 \big\rangle_T &= \frac{2}{d^2 \pi^2}
\sum\limits_{\substack{\vec{k} \in \Lambda ^ {*} \setminus \{ 0 \}
\\ |\vec{k}| \le T^{3/2}}} \frac{\sin ^2 \bigg( \frac{\pi
|\vec{k}|}{L} \bigg) } {|\vec{k}| ^{3}}
\bigg(1- \hat{\psi}  \big( \frac{|\vec{k}|}{\sqrt {M}} \big) \bigg) ^ 2 \\
&+O \Bigg( \sum\limits_{\substack{\vec{k},\, \vec{l} \in \Lambda ^ {*} \setminus \{ 0 \} \\
|\vec{k}| \ne |\vec{l}| \le T^{3/2}}} \frac{ 1} {|\vec{k}| ^{3/2}
|\vec{l}|^{3/2}} \hat{\omega} \bigg( T \big( |\vec{k}| - |\vec{l}|
\big) \bigg) \Bigg)
\end{split}
\end{equation}

We will evaluate the diagonal contribution now. For $|\vec{k}| \le
\sqrt{M}$,
\begin{equation*}
\hat{\psi}  \bigg( \frac{|\vec{k}|}{\sqrt {M}} \bigg) = 1 + O
\bigg(  \frac{|\vec{k}|}{\sqrt {M}} \bigg),
\end{equation*}
and so the diagonal contribution is:
\begin{equation*}
\frac{1}{M} \sum\limits_{\substack{\vec{k} \in \Lambda ^ {*}
\setminus \{ 0 \} \\ 1 \ll |\vec{k}| \le \sqrt{M}}}
\frac{1}{|\vec{k}|} + \sum\limits_{\substack{\vec{k} \in \Lambda ^
{*} \setminus \{ 0 \} \\ \sqrt{M} \le |\vec{k}| \le T^{3/2}}}
\frac{1 } {|\vec{k}| ^{3}} \ll \frac{1}{\sqrt{M}},
\end{equation*}
converting the sums into corresponding integrals and evaluating
these integrals in the elliptic variables.

Finally, we are evaluating the off-diagonal contribution to
\eqref{eq:asmp eqlty P} (that is, the second sum in the right-hand
side of \eqref{eq:asmp eqlty P}). Set $0 < \delta_0 < 1$. With no
loss of generality, we may assume that $|\vec{k}| < |\vec{l}|$.
Evaluating the contribution of pairs $\vec{k},\, \vec{l}$ with
\begin{equation*}
|\vec{l}| ^2 - |\vec{k}| ^2 \ge \frac{|\vec{k}|}{T^{1-\delta_0}}
\end{equation*}
gives:

\begin{equation*}
\ll \sum\limits_{\substack{\vec{k},\, \vec{l} \in \Lambda ^ {*} \setminus \{ 0 \} \\
|\vec{k}| < |\vec{l}| \le T^{3/2}}} \frac{ 1} {|\vec{k}| ^{3/2}
|\vec{l}|^{3/2}} \hat{\omega} \bigg( T \big( |\vec{k}| - |\vec{l}|
\big) \bigg) \ll T^{-A\delta_0 + 6}
\end{equation*}
for every $A>0$, since
\begin{equation*}
T \big( |\vec{l}| - |\vec{k}| \big) = T \frac {|\vec{l}| ^2 -
|\vec{k}| ^2}{  |\vec{k}| + |\vec{l}|  } \gg T^{\delta_0}
\frac{|\vec{k}|}{|\vec{k}| + |\vec{l}| } \ge T^{\delta_0}
\frac{|\vec{k}|}{|\vec{k}| + 2|\vec{k}| } \gg T^{\delta_0},
\end{equation*}
as otherwise,
\begin{equation*}
T \big( |\vec{l}| - |\vec{k}| \big) \ge T \big( |\vec{k}| ) \gg T
\gg T^{\delta_0}.
\end{equation*}
Thus the contribution of such terms is negligible.

In order to bound the contribution of pairs of
$\Lambda^{*}$-vectors with
\begin{equation*}
|\vec{l}| ^2 - |\vec{k}| ^2 \le \frac{|\vec{k}|}{T^{1-\delta_0}}
\end{equation*}
we use the Diophantinity of $\beta$ again. Recall that we chose
$\delta_{\Lambda ^* } (y) = \frac{c_2}{y^{K_0}}$ with a constant
$c_2$ in the beginning of the current proof. Choose a constant
$R_0 > 0$ and assume that $|\vec{l}|^2 \le c L^{R_0}$, for a
constant $c$. Then
\begin{equation*}
\begin{split}
|\vec{l}| ^2 - |\vec{k}| ^2 &\ge \delta_{\Lambda^*} (L^{R_0}) \gg
\frac{1}{L^{K_0 R_0}} \gg \frac{1}{M^{K_0 R_0 / 2}} \gg
\frac{|\vec{k}|}{T^{1-\delta_0}}.
\end{split}
\end{equation*}
Therefore, for an appropriate choice of $c$, there are no such
pairs. Denote
\begin{equation*}
S_n :=\bigg\{ (\vec{k},\, \vec{l})\in (\Lambda^*) ^ 2 :\: 2^n \le
|\vec{k}|^2 \le 2^{n+1},\, |\vec{k}|^2 \le |\vec{l}|^2 \le
|\vec{k}|^2 + \frac{2^{n/2}}{T^{1-\delta_0}} \bigg\}
\end{equation*}
Thus, by dyadic partition, the contribution is:
\begin{equation*}
\begin{split}
&\ll\sum\limits_{n = \lfloor R_0 \log L \rfloor}^{\lceil 3 \log T
\rceil }
\sum\limits_{\substack{ 2^n \le |\vec{k}|^2 \le 2^{n+1} \\
|\vec{k}|^2 \le |\vec{l}|^2 \le |\vec{k}|^2 +
\frac{2^{n/2}}{T^{1-\delta_0}}} } \frac{1}{|\vec{k}|^{3/2}
|\vec{l}|^{3/2}} \hat{\omega} \bigg( T \big( |\vec{k}| - |\vec{l}|
\big) \bigg)  \\ &\ll
 \sum\limits_{n = \lfloor R_0 \log L \rfloor}^{\lceil 3 \log T \rceil}
\frac{ \# S_n }{2^{3n/2}},
\end{split}
\end{equation*}
using $|\hat{\omega}| \ll 1$ everywhere. 
In order to bound the size of $S_n$, we use the following lemma,
which is just a restatement of lemma 3.1 from ~\cite{BL}. We will
prove it immediately after we finish proving proposition
\ref{prop:bnd 2nd mom diff}.
\begin{lemma}
\label{lem:S_n bnd blh} Let $\Lambda = \big \langle 1,\, i\eta
\big\rangle$ be a rectangular lattice. Denote
\begin{equation*}
A(R,\delta) := \{ (\vec{k},\, \vec{l}) \in \Lambda :\: R \le
|\vec{k}|^2 \le 2R,\, |\vec{k}|^2 \le |\vec{l}|^2 \le |\vec{k}|^2
+ \delta \}.
\end{equation*}
Then if $\delta > 1$, we have for every $\epsilon>0$,
\begin{equation*}
\# A(R,\delta) \ll _{\epsilon} R^{\epsilon} \cdot R \delta
\end{equation*}
\end{lemma}


Thus, lemma \ref{lem:S_n bnd blh} implies
\begin{equation*}
\# S_n \\ \ll 2^{n+\epsilon (n/2)} \max \bigg( 1, \,
\frac{2^{n/2}}{T^{1-\delta_0}} \bigg),
\end{equation*}
for every $\epsilon > 0$. Thus the contribution is:
\begin{equation*}
\begin{split}
\ll \sum\limits_{n= R_0 \log L - 1} ^{C \log {T} + 1}
\frac{1}{2^{3n/2}} \cdot 2^{n+\epsilon n/2} \cdot 1 &+
\sum\limits_{n= C \log T - 1} ^{3 \log {T} + 1} \frac{1}{2^{3n/2}}
\cdot 2^{n+\epsilon n/2}
\cdot \frac{2^{n/2}}{T^{1-\delta_0}} \\
&\ll L^{-R_0 (1-\epsilon)/2 } + \frac{\log T}{T^{1-\delta_1}} \ll
L^{-R_0 (1-\epsilon) /2 },
\end{split}
\end{equation*}
since $L$ is much smaller than $T$. Since $R_0$ is arbitrary, and
we have assumed $M = O (L^{s_0})$, that implies
\begin{equation*}
\big\langle P_{\Lambda} ^2 \big\rangle_T \ll \frac{1}{\sqrt{M}}.
\end{equation*}

Collecting all our results, and using them on \eqref{eq:CS ineq
for var} we obtain
\begin{equation*}
\Bigg\langle \bigg| S_{\Lambda} (t,\, \rho ) -
\tilde{S}_{\Lambda, \, M, \, L} (t) \bigg| ^ {2}  \Bigg\rangle \ll
\frac{1}{\sqrt{M}} + \frac{1}{T^{1-\delta'}} +
\frac{\sqrt{\log{M}}}{M^{1/4} T^{1/2 - \delta' / 2}} \ll
\frac{1}{\sqrt{M}},
\end{equation*}
again, since $M$ is much smaller than $T$.
\end{proof}

\begin{proof}[Proof of lemma \ref{lem:S_n bnd blh}]
Let $\vec{k} = (k_1,\, i\eta k_2)$ and $\vec{l} = (l_1,\, i\eta
l_2)$. Denote $\mu:=\eta^2$, $n:=l_1^2-k_1^2$ and
$m:=k_2^2-l_2^2$. The number of 4-tuples $(k_1,\, k_2,\, l_1,\,
l_2)$ with $m \ne 0$ is
\begin{equation*}
\# A(\delta,\, T) \ll \sum\limits_{\substack{0 \le n-\mu m \le
\delta \\ 1 \le m \le 4R}} d(n)d(m) \ll \delta \sum\limits_{1 \le
m \le 4R} d(m)^2 \ll R^{1+\epsilon} \delta
\end{equation*}
Next, we bound the number of 4-tuples with $m=0,\, n\ne 0$:
\begin{equation*}
\sum\limits_{k_2=0}^{\sqrt{2R}} \sum\limits_{0 < n <\delta} d(n)
\ll R^{1/2+\epsilon} \delta,
\end{equation*}
and similarly we bound the number of 4-tuples with $n=0,\, m\ne
0$.

All in all, we have proved that
\begin{equation*}
\# A(\delta,\, T) \ll R^{1+\epsilon} \delta
\end{equation*}
\end{proof}

From now on we will assume that $\Lambda = \big\langle 1,\,
i\alpha \big\rangle$ with a Diophantine $\gamma := \alpha^2$, and
so the use of proposition \ref{prop:bnd 2nd mom diff} is
justified.

\begin{lemma}
\label{lem:pr big diff shrp smth} Under the conditions of
proposition \ref{prop:bnd 2nd mom diff}, for all fixed $\xi > 0$,
\begin{equation*}
\PP _{\omega,\, T} \Bigg\{ \bigg|\frac{S_{\Lambda} (t,\, \rho )
}{\sigma} - \frac{\tilde{S}_{\Lambda,\, M,\, L } (t)}{\sigma}
\bigg|  > \xi \Bigg\} \rightarrow 0,
\end{equation*}
as $T \rightarrow \infty$, where $\sigma^2 = \frac{8 \pi }{d L}$.
\end{lemma}
\begin{proof}
Use Chebychev's inequality and proposition \ref{prop:bnd 2nd mom
diff}.
\end{proof}

\begin{corollary}
\label{cor:main res omg} For a number $\alpha\in\R$, suppose that
$\alpha^2$ is \underline{strongly Diophantine} and denote $\Lambda
= \big\langle 1,\, \alpha \big\rangle$. Then if $L \rightarrow
\infty$, but $L=O \big( T^{\delta} \big)$ for all $\delta > 0$ as
$T\rightarrow \infty$, then for any interval $\A$,
\begin{equation*}
\PP _{\omega,\, T} \bigg\{ \frac{S_{\Lambda} (t,\, \rho )
}{\sigma} \in \A \bigg\} \rightarrow \frac{1}{\sqrt{2\pi }}
\int\limits_{A} e ^ {-\frac{x^2}{2}} dx,
\end{equation*}
where $\sigma^2 = \frac{8 \pi }{d L}$.
\end{corollary}
\begin{proof}
Set $M=L^3$, then, obviously, $L,\, M$ satisfy the conditions of
lemma \ref{lem:pr big diff shrp smth} and theorem \ref{thm:norm
dist smth}. Denote $X(t) := \frac{S_{\Lambda} (t,\, \rho
)}{\sigma}$ and $Y(t) := \frac{\tilde{S}_{\Lambda,\, M }
(t)}{\sigma}$. In the new notations lemma \ref{lem:pr big diff
shrp smth} states that for any $\xi > 0 $,
\begin{equation}
\label{eq:prob diff X Y 0} \PP _{\omega,\, T} \big\{ |X(t) - Y(t)|
> \xi \big\} \rightarrow 0,
\end{equation}
as $T \rightarrow \infty$. Now, for every $\epsilon > 0$,
\begin{equation*}
\big\{a \le X \le b \big\} \subseteq \big\{ a- \epsilon \le Y \le
b + \epsilon \big\} \cup \big\{ |X-Y| > \epsilon \big\},
\end{equation*}
and so, taking $\limsup\limits_{T\rightarrow\infty}
{\PP_{\omega,\, T}}$ of both of the sides, we obtain:
\begin{equation*}
\begin{split}
\limsup\limits_{T\rightarrow\infty} {\PP_{\omega,\, T} \big\{a \le
X \le b \big\}} &\le \lim\limits_{T\rightarrow\infty}
{\PP_{\omega,\, T}\big\{ a- \epsilon \le Y \le b + \epsilon
\big\}} \\ &= \frac{1}{\sqrt{2\pi }}
\int\limits_{a-\epsilon}^{b+\epsilon} e ^ {-\frac{x^2}{2}} dx,
\end{split}
\end{equation*}
due to \eqref{eq:prob diff X Y 0} and theorem \ref{thm:norm dist
smth}. Starting from
\begin{equation*}
\big\{a + \epsilon \le Y \le b - \epsilon \big\} \subseteq \big\{
a \le X \le b \big\} \cup \big\{ |X-Y| > \epsilon \big\},
\end{equation*}
and doing the same manipulations as before, we get the converse
inequality, and thus this implies the result of the present
corollary.

\end{proof}

We are now in a position to prove our main result, namely, theorem
\ref{thm:norm dist}. It states that the result of corollary
\ref{cor:main res omg} holds for $\omega = \1_{[1,\,2]}$, the
indicator function. We are unable to substitute it directly
because of the rapid decay assumption on $\hat{\omega}$.
Nonetheless, we are able to prove the validity of the result by
the means of {\em approximating} the indicator function with
functions which will obey the rapid decay assumption. The proof is
essentially the same as of theorem 1.1 in ~\cite{HR}, pages
655-656, and we repeat it in this paper for the sake of the
completeness.

\begin{proof}[Proof of theorem \ref{thm:norm dist}]
Fix $\epsilon > 0$ and approximate the indicator function
$\1_{[1,\,2]}$ above and below by smooth functions $\chi_{\pm} \ge
0$ so that $\chi_{-} \le \1_{[1,\,2]} \le \chi_{+}$, where both
$\chi_{\pm}$ and their Fourier transforms are smooth and of rapid
decay, and so that their total masses are within $\epsilon$ of
unity $\big| \int \chi_{\pm} (x) dx - 1 \big| < \epsilon $. Now,
set $\omega_{\pm} := \chi_{\pm}/ \int\chi_{\pm}$. Then
$\omega_{\pm}$ are "admissible", and for all $t$,
\begin{equation}
\label{eq:om approx qlty} (1-\epsilon ) \omega_{-} (t) \le
\1_{[1,\,2]} (t) \le (1+\epsilon ) \omega_{+} (t).
\end{equation}
Now,
\begin{equation*}
meas \bigg\{ t\in [T,\, 2T ] :\: \frac{S_\Lambda (t,\, \rho ) }
{\sigma } \in \A \bigg\} = \int\limits_{-\infty}^{\infty} \1_{\A}
\bigg( \frac{S_\Lambda (t,\, \rho ) } {\sigma } \bigg)
\1_{[1,\,2]} \big( \frac{t}{T} \big) dt,
\end{equation*}
and since \eqref{eq:om approx qlty} holds, we find that
\begin{equation*}
\begin{split}
(1-\epsilon) \PP_{\omega_{-},\, T} \bigg\{ \frac{S_{\Lambda, \, M,
\, L} } {\sigma} \in \A \bigg\}
 &\le \frac{1}{T} meas \bigg\{ t\in [T,\, 2T ] :\: \frac{S_\Lambda (t,\, \rho ) } {\sigma } \in \A \bigg\} \\
&\le (1+\epsilon) \PP_{\omega_{+} ,\, T} \bigg\{ \frac{S_{\Lambda,
\, M, \, L} } {\sigma } \in \A \bigg\}.
\end{split}
\end{equation*}
As it was mentioned immediately after the definition of the strong
Diophantinity property, $\alpha$'s being strongly Diophantine
implies the same for $\alpha^2$, making a use of corollary
\ref{cor:main res omg} legitimate. Now by corollary \ref{cor:main
res omg}, the two extreme sides of the last inequality have a
limit, as $T\rightarrow\infty$, of
\begin{equation*}
(1\pm\epsilon) \frac{1}{\sqrt{2\pi }} \int\limits_{A} e ^
{-\frac{x^2}{2}} dx,
\end{equation*}
and so we get that
\begin{equation*}
(1-\epsilon) \int\limits_{A} e ^ {-\frac{x^2}{2}} dx \le
\liminf\limits_{T\rightarrow\infty} \frac{1}{T} meas \bigg\{ t\in
[T,\, 2T ] :\: \frac{S_\Lambda (t,\, \rho ) } {\sigma } \in \A
\bigg\}
\end{equation*}
with a similar statement for $\limsup$; since $\epsilon > 0$ is
arbitrary, this shows that the limit exists and equals
\begin{equation*}
\lim\limits_{T\rightarrow\infty} \frac{1}{T} meas \bigg\{ t\in
[T,\, 2T ] :\: \frac{S_\Lambda (t,\, \rho ) } {\sigma } \in \A
\bigg\} = \frac{1}{\sqrt{2\pi }} \int\limits_{A} e ^
{-\frac{x^2}{2}} dx,
\end{equation*}
which is the Gaussian law.
\end{proof}

\paragraph{Acknowledgement.}
This work was supported in part by the EC TMR network
\textit{Mathematical Aspects of Quantum Chaos}, EC-contract no
HPRN-CT-2000-00103 and the Israel Science Foundation founded by
the Israel Academy of Sciences and Humanities. This work was
carried out as part of the author's PHD thesis at Tel Aviv
University, under the supervision of prof. Ze\'{e}v Rudnick. The
author wishes to thank the referee for many useful comments. A
substantial part of this work was done during the author's visit
to the university of Bristol.

\end{document}